\documentclass[reqno,12pt,a4paper]{amsart}

\usepackage{mathrsfs}
\usepackage{amssymb}
\usepackage{amsfonts}
\usepackage{latexsym}
\usepackage{amsthm}

\usepackage{graphicx}
\def\lb{\label}

\newcommand{\er}[1]{\textrm{(\ref{#1})}}

\theoremstyle{plain}
\newtheorem{theorem}{\bf Theorem}[section]
\newtheorem{lemma}[theorem]{\bf Lemma}
\newtheorem{corollary}[theorem]{\bf Corollary}
\newtheorem{proposition}[theorem]{\bf Proposition}
\newtheorem{definition}[theorem]{\bf Definition}

         \def\mA{{\mathscr A}}
   \def\cB{{\mathcal B}}       \def\mB{{\mathscr B}}
\def\g{\gamma}         
\def\G{\Gamma}  \def\cD{{\mathcal D}}       \def\mD{{\mathscr D}}
\def\d{\delta}  \def\cE{{\mathcal E}}       
  \def\cF{{\mathcal F}}       
\def\c{\chi}           
\def\z{\zeta}   \def\cH{{\mathcal H}}       \def\mH{{\mathscr H}}
\def\e{\eta}           
    \def\cJ{{\mathcal J}}       
 \def\cK{{\mathcal K}}       
  \def\cL{{\mathcal L}}       
\def\l{\lambda} \def\cM{{\mathcal M}}       \def\mM{{\mathscr M}}
\def\L{\Lambda} \def\cN{{\mathcal N}}       \def\mN{{\mathscr N}}
     \def\cO{{\mathcal O}}

\def\s{\sigma}         \def\mR{{\mathscr R}}
  \def\cS{{\mathcal S}}       \def\mS{{\mathscr S}}
           
\def\f{\phi}           
\def\F{\Phi}           
           \def\mW{{\mathscr W}}
\def\o{\omega}  \def\cX{{\mathcal X}}       
            
     \def\cZ{{\mathcal Z}}       
\def\O{\Omega}

\newcommand{\gF}{\mathfrak{F}}           \newcommand{\gf}{\mathfrak{f}}

\def\ve{\varepsilon}   \def\vt{\vartheta}    \def\vp{\varphi}    \def\vk{\varkappa}

  \def\C{{\mathbb C}}
        \def\N{{\mathbb N}}   \def\R{{\mathbb R}}      \def\Z{{\mathbb Z}}

\def\lt{\biggl}                  \def\rt{\biggr}
\def\ol{\overline}               \def\wt{\widetilde}
\def\no{\noindent}

\let\ge\geqslant                 \let\le\leqslant
                
\def\/{\over}                    \def\iy{\infty}
\def\sm{\setminus}               
\def\ss{\subset}                 \def\ts{\times}
\def\pa{\partial}                \def\os{\oplus}

\def\el2{\ell^{\,2}}             \def\1{1\!\!1}

\def\co{zj_{\vk-1}(z)}
\def\si{zj_{\vk}(z)}

\def\arg{\mathop{\mathrm{arg}}\nolimits}
\def\const{\mathop{\mathrm{const}}\nolimits}
\def\det{\mathop{\mathrm{det}}\nolimits}

\def\Im{\mathop{\mathrm{Im}}\nolimits}

\def\Re{\mathop{\mathrm{Re}}\nolimits}

\def\sign{\mathop{\mathrm{sign}}\nolimits}

\def\supp{\mathop{\mathrm{supp}}\nolimits}
\def\Tr{\mathop{\mathrm{Tr}}\nolimits}
\def\BBox{\hspace{1mm}\vrule height6pt width5.5pt depth0pt \hspace{6pt}}

\def\qqq{\qquad}
\def\qq{\quad}
\let\ge\geqslant
\let\le\leqslant
\let\geq\geqslant
\let\leq\leqslant
\newcommand{\ca}{\begin{cases}}
\newcommand{\ac}{\end{cases}}
\newcommand{\ma}{\begin{pmatrix}}
\newcommand{\am}{\end{pmatrix}}
\renewcommand{\[}{\begin{equation}}
\renewcommand{\]}{\end{equation}}
\def\eq{\begin{equation}}
\def\qe{\end{equation}}
\def\[{\begin{equation}}

\begin{document}
\bibliographystyle{plain}
%{alpha}

\title[{Resonances for  the radial Dirac operators}]
{Resonances for  the radial Dirac operators}

\date{\today}

\author[ Alexei Iantchenko]{ Alexei Iantchenko}
\address{Malm{\"o} h{\"o}gskola, Teknik och samh{\"a}lle, 205 06
  Malm{\"o}, Sweden, email: ai@mah.se }
\author[Evgeny Korotyaev]{Evgeny Korotyaev}
\address{Mathematical Physics Department, Faculty of Physics, Ulianovskaya 2,
St. Petersburg State University, St. Petersburg, 198904, Russia,
 \ korotyaev@gmail.com.
}


\begin{abstract}
\no We consider the radial  Dirac operator     with compactly supported potentials.
We study  resonances as the poles of scattering matrix or equivalently as the zeros of modified Fredholm determinant.
We obtain the following properties of the resonances:
1) asymptotics of counting function,
2)  in the massless case we get the trace formula in terms of resonances.
\end{abstract}

\maketitle

%%%%%%%%%%%%%%%%%%%%%%%%%%%%%%%%%%%%%%%%%%%%
%% MAINMATTER
%%%%%%%%%%%%%%%%%%%%%%%%%%%%%%%%%%%%%%%%%%%%

\vskip 0.25cm

\section{Introduction}

 The spherically symmetric Dirac operator in $\R^3$ (in the units $\hbar=c=1$) has partial-wave decomposition in 1D radial Dirac operators
\[
\lb{Dirac3D}
-i\sum_{j=1}^3\alpha\cdot\nabla +\beta m+V(|\overline{x}|)\cong\bigoplus_{\vk\in\Z\setminus\{0\}}\,\,
\bigoplus_{m_\vk={1\/2}-|\vk|}^{|\vk|-\frac12}\left(-i\s_2\frac{d}{dr}+\sigma_3 m+\sigma_1\frac{\vk}{r}+V(r)\right),
\] 
where
 $\alpha=(\alpha_1,\alpha_2,\alpha_3),$ and
$\alpha_j,\beta$ are the $4\times 4$ Dirac matrices
$$
\alpha_j=\ma 0 &\s_j\\ \s_j &0\am,\qq \beta=\ma I_2 &0\\ 0 &-I_2\am,\qqq
j=1,2,3
$$ 
and $2\times 2$ matrices $\s_1, \s_2, \s_3$ are the Pauli matrices given by
$$
\sigma_1=\ma 0 &1\\1&0 \am,\qq \sigma_2=\ma 0 &-i\\i&0 \am,\qq\s_3=\ma 1 & 0\\0&-1 \am.
$$
We put also $\sigma_0=I_2$ the $2\ts 2$ identity matrix and $m>0$ is the mass, $V(|\overline{x}|)=v(r)\sigma_0\in L^1(\R_+)$ is spherically symmetric electrostatic field,
% and $\supp v \subset [0,\gamma],$
$\vk$ is the spin-orbit coupling parameter satisfying \[\lb{defvk}\vk=\pm (j+\frac12)\,\,\mbox{ if }\,\, \ell=j\pm\frac12,\] where $j=\frac12,\frac32,\frac52,\ldots$ and $\ell=j\pm\frac12$ are the total and orbital angular momentum numbers respectively.
Relation (\ref{defvk}) is usually taken as definition of $\vk$ (see \cite{BjorkenDrell1964}, \cite{GriesemerLutgen1999}), it says that the sign of $\vk$ indicates whether spin and orbital angular momentum of the upper component are ``parallel'' or ``anti-parallel.''

As the spectral characteristics of the 1D radial Dirac operator
$H =-i\s_2\pa_x+\sigma_3 m+\sigma_1\frac{\vk}{r}+v(r)\sigma_0$
only depend on $|\vk |,$ so it is enough to suppose that $\vk >0.$ 
%Then the multiplicity  is $2\cdot 2\vk.$

In this paper we will study the scattering resonances. Resonances are the complex numbers associated to  the outgoing modes and  can be defined as the poles of analytic continuation of the resolvent acting between suitable distribution spaces of distributions.  (see Definition \ref{DefRes} below). 
From a physicists
point of view, the resonances were first studied by Regge in 1958 (see \cite{Regge1958}). Since then,
 the properties of  resonances for the Schr\"odinger type operators had been the object of intense study and we refer to    \cite{SjostrandZworski1991} and \cite{Zworski1994} for the mathematical approach in the multi-dimensional case  and references
given there.The resonances were defined by the method   of complex scaling under the hypothesis that a real-valued smooth potential extends analytically to a complex conic neighborhood of the real domain at infinity and tends to $0$ sufficiently fast there as $x\rightarrow\infty.$   As result, only local or semi-classical properties of resonances could be derived.  In the multi-dimensional Dirac case resonances were studied locally in \cite{BalslevHelffer1992}.

We are interested in the global properties of resonances which imposes further restrictions on the potential. The potential is supposed to have compact support or,  at least, super-exponentially decreasing at infinity.  In this context,  the resonances for the 1D Schr\"odinger operator are well studied, see
Froese \cite{Froese1997}, Simon \cite{Simon2000}, Korotyaev \cite{Korotyaev2004},  Zworski
\cite{Zworski1987} and references given there. We recall that Zworski \cite{Zworski1987}
obtained the first results about the asymptotic distribution of
resonances for the Schr\"odinger operator with compactly supported
potentials on the real line. Different properties of resonances were
determined in \cite{Hitrik1999} and \cite{Korotyaev2011}. Inverse problems (characterization,
recovering, plus uniqueness) in terms of resonances were solved by
Korotyaev  for the Schr\"odinger operator with a compactly supported
potential on the real line \cite{Korotyaev2005} and the half-line \cite{Korotyaev2004}.
The "local resonance" stability problems were considered in
\cite{Korotyaev2004a}, \cite{Marlettaetal2010}.

 Similar questions for Dirac operators are much less
studied. However, there are a number of papers dealing with other related problems (see \cite{IantchenkoKorotyaev2013a} for the references). 

 In
\cite{IantchenkoKorotyaev2013} we consider the 1D massless Dirac operator    on the
real line with compactly supported potentials. It is a special kind
of the Zakharov-Shabat operator (see \cite{DEGM}, \cite{Novikovetal1984}). 
Technically, this case is simpler than the massive Dirac
operator, since in the massless case the Riemann
surface consists of two disjoint sheets $\C.$   Moreover, the resolvent
has a simple representation.  The goal of  \cite{IantchenkoKorotyaev2013} was to give a clear untechnical presentation of ideas which are generalized  \cite{IantchenkoKorotyaev2013a} and in the present paper and will be further developed in our other papers  in preparation  \cite{IantchenkoKorotyaev2014b},  \cite{IantchenkoKorotyaev2014c}.
Moreover, in \cite{IantchenkoKorotyaev2013} we were even able  to
prove the trace formulas in terms of resonances. Similar results are obtained in the present paper for $m=0$ (see Theorem \ref{T4} ). We have not been able to get  a similar result in the general situation with non-zero mass. Note that in the
massless case the relation between the modified Fredholm determinant
$D$ and the Jost function $f_1^+(0,\l)$ (corresponding to $a$ for
the problem on the line in \cite{IantchenkoKorotyaev2013}, the inverse of the
transmission coefficient) is much easier than in the massive case
(see Theorem \ref{T1}), namely $D(\l)=a(\l),$ with no 
proportionality factor in between. Note that in the singular case as discussed in the present paper this is no longer true, even in the massless case (see Theorem \ref{T1})

In \cite{IantchenkoKorotyaev2013a} we consider the regular case which corresponds to radial Dirac operator $H$ without singular potential $\vk/r$ (i.e. $\vk=0$) and general perturbation potential $$
V(x)=\ma p_1 & q\\ q & p_2\am (x),\qq x\geq 0
$$
with real-valued functions $p_1,$ $p_2$ and $q.$ The present paper concerns  the singular at $x=0$ problem,  $\vk\neq 0.$ In comparison to   \cite{IantchenkoKorotyaev2013},  the techniques used in the present paper are heavier due to the use of Bessel functions, and the asymptotics are more complicated due to the presence of several (small or large) parameters.

\section{Definitions and main results}

\subsection{Modified Fredholm determinant.}
We will write $x$ instead of $r=|\overline{x}|.$
We consider the radial  Dirac operator $H=H_0+V$ acting on the Hilbert space $L^2(\R_+ )\os L^2(\R_+ )$, where  $H_0$ is the free radial Dirac operator 
 given by 
 \[
\lb{free_radial}
H_0 f=\rt(-i\s_2\pa_x+\sigma_3 m+\sigma_1\frac{\vk}{x}\rt)f=\ma m & -\partial_x+\frac{\vk}{x}\\\partial_x+\frac{\vk}{x} & -m\am f,\qqq
f=\ma f_1\\ f_2\am,
\] 
where  $\vk\in\Z_+=\{1,2,3,\ldots\}$ and  $f$ satisfies the Dirichlet boundary condition
\[
\lb{Dirichlet_cond} f_1(0)=0.
\] 
Here $V$ is the real diagonal  matrix-valued potential, satisfying the following conditions:
\[
\lb{weakas} 
 V=v I_2,\qqq    \int_0^\infty(1+x)|v(x)|dx <\infty.
\]

The boundary condition  \er{Dirichlet_cond} and our  assumption \er{weakas} on $V$, guaranty that the differential operator $H$ is  self-adjoint on the Hilbert space $L^2(\R_+ )\os L^2(\R_+ )$.
 The spectrum of $H_0$  is absolutely continuous and is given by
$$
\s(H_0)=\s_{\rm ac}(H_0)=\R\sm (-m,m).
$$ 
The spectrum of $H$ consists of the absolutely continuous part $\s_{\rm ac}(H)=\s_{\rm ac}(H_0)$ plus a finite number of simple eigenvalues in the gap $(-m,m).$

It is well known that the wave operators $W_\pm=W_\pm(H,H_0)$ for
the pair $H_0, H$ given by
$$
W_\pm=s\!-\!\lim e^{itH}e^{-itH_0} \qqq \mbox{as} \qqq t\to \pm\iy,
$$
exist and  are complete (even under much less restrictive assumptions on the potential than considered here, see \cite{Thaller1992}).
Thus the scattering operator
$S=W_+^*W_-$ is unitary. The operators $H_0$ and $S$ commute and thus
are simultaneously diagonalizable:
\[
\lb{DLH0}
L^2(\R_+)\os L^2(\R_+)=\int_\R^\oplus \mH_\l d\l,\qqq H_0=\int_\R^\oplus\l I_\l d\l,\qqq
S=\int_\R^\oplus \cS(\l)d\l;
\]
here $I_\l$ is the identity in the fiber space $\mH_\l=\C$ and $\cS(\l)$ is
the scattering matrix (which is a  scalar function of $\l\in \R$ in our case)
 for the pair $H_0, H$ (see \cite{Thaller1992}).

Now, we introduce a basis of {\bf Jost solutions}  $f^\pm$ for $H$ by the conditions
  \[\lb{JostAsy1-0}
\begin{aligned}
& H  f^\pm=\l f^\pm,\qq
  f^\pm(x,\vk,\l)= (\mp i k)^\vk e^{\pm i k(\l) x}\ma\pm k_0(\l)\\ 1\am+o(1),\,\,\mbox{as}\,\,x\rightarrow\infty,\\ & k_0(\l)=\frac{\l +m}{i k(\l)},\qq k(\l)=\sqrt{\l^2-m^2},\qq \l\in\sigma_{\rm ac}(H_\vk),\end{aligned}\]
 where the function $k(\l)$ is quasi-momentum and defined below in (\ref{q-m}). Note that $f^-(x,\l)=\overline{f^+(x,\l)}$ for $\l\in\sigma_{\rm ac}(H).$
The Jost solutions for the unperturbed system ($v=0,$ associated with free radial Dirac operator (\ref{free_radial})) are defined by the same condition (\ref{JostAsy1}) and are denoted by $\psi^\pm(x,\l).$
The {\bf Jost function} is given by
\[\lb{Jost_function_def}\gf^+(\l)=\lim_{x\rightarrow 0}\frac{x^\vk}{(2\vk-1)!!}f_1^+(x,\l).\]
We denote $\gf^{0,+}(\l)$  the Jost function for 
 the unperturbed Dirac system ($v=0$). We show that  $\gf^{0,+}(\l)=k_0(\l).$

From  results  in \cite{Blancarteetal1995} recalled  in Theorem \ref{th-BFW2} it follows that (under appropriate conditions on $v$) \[\lb{Jost_inf}\gf^+(i\eta)\rightarrow e^{i\left(\int_0^\infty v(t)dt-\frac{\pi}{2}\right)},\qq\mbox{as}\qq \eta\rightarrow \infty,\] and therefore we take the unique branch $\log \left(e^{i\left(\frac{\pi}{2}-\int_0^\infty v(t)dt\right)}\gf^+(\l)\right)=o(1)$ as $\l=i\eta,$ $\eta\rightarrow\infty.$
Due to (\ref{Jost_inf}) we can define  the unique branch $\log
\gf^+(\l)$ in $\C_+$ and define the functions
$$
\log \gf^+(\l)=\log |\gf^+(\l )|+i\arg \gf^+(\l ),\qq \l\in \C_+ ,
$$
where the function $\f_{\rm sc}=\arg \gf^+(\l )+\pi/2$ is called the
scattering phase (or the spectral shift function).

The scattering matrix $\cS(\l),$ $\l\in\s_{\rm ac}(H_0),$ for the pair $H,H_0$   is then given by
$$\cS(\l)=-\frac{\overline{\gf^+(\l+i0)}}{\gf^+(\l+i0)}=e^{-2i\phi_{\rm sc}},\qq \mbox{for}\,\,\l\in\s_{\rm ac}(H_0).$$

The minus sign comes from our choice of the normalization of the Jost
solutions at the spatial infinity  (\ref{JostAsy1}).
 Property (\ref{Jost_inf}) implies
$$
\phi_{\rm sc}(\l)=\int_0^\infty v(t)dt+o(1),\qq\mbox{as}\qq \Im\l\rightarrow\infty.
$$

 We will show below that the Jost function and scattering matrix is related to the
 {\em modified Fredholm determinant}  introduced as follows. We
set
\[
\lb{Y0}
\begin{aligned}
&V=|V|^{1\/2}V^{1\/2},\qqq V^{1\/2}=|v|^{1\/2}I_2 \sign v,
\\
& R_0(\l)=(H_0-\l)^{-1},\qqq  Y_0(\l)=|V|^{1\/2}\, R_0(\l)\,
V^{1\/2}, \qqq \l \in {\C}_\pm.
\end{aligned}
\]
Here $\C_\pm=\{\l\in\C: \pm\Im \l>0\}$ denote the upper and lower
half plane and $\l$ is a spectral parameter. 
Observing  that the operator valued function $Y_0(\l)$ is in the Hilbert-Schmidt class $\cB_2$ but
 not in the trace class $\cB_1,$ (see \cite{Korotyaev2014}), we define the modified Fredholm determinant $D(\l)$ (see \cite{GohbergKrein1969}) by
$$
D(\l)=\det\left[ (I+Y_0(\l)) e^{-Y_0(\l)}\right], \qqq \forall\,\l\in \C_+.
$$
We will show later that the function   
\[\lb{Omega}
\Omega(\l)=\Tr (Y_0(\l+i0)-Y_0(\l-i0))\qq
\mbox{if}\qq \l\in \R\sm \{\pm m\}
\] 
is well defined. Note that $\O(\l)=0$ on the interval $(-m,m)$.

We formulate the main results of this paper connecting the modified Fredholm determinant $D$ and the Jost function $\gf^+.$

\begin{theorem}\lb{T1}
Let $v\in L^1(\R_+)\cap L^2(\R_+)$ and $v\in L^\infty(0,a)$ for some $a>0.$ Then  the Jost function $\gf^+(\l)$ and the
determinant $D(\l)$ are analytic in $\C_+,$  continuous up to
$\R\sm \{\pm m\}$  and satisfy
 \[
\lb{Df}
\begin{aligned}
&\cS(\l)=\frac{D(\l-i 0)}{D(\l+i 0)}\,e^{-2i\Omega(\l)}=e^{-2i\f_{\rm sc}(\l)},\\
&\f_{\rm sc}(\l)=\Omega(\l)+\arg D(\l+i0),\qq\forall\l\in\s_{\rm ac}(H_0),\qq
\l\neq\pm m.
\end{aligned}
\]
Here the function $\O$ (defined in (\ref{Omega})) is continuous on
$\R\sm \{\pm m\}$ and satisfies
\[
\lb{Oid}\Omega(\l)=\int_0^\infty v(y)\left( \frac{k}{\l-m}\left[kyj_\vk(ky)\right]^2 +\frac{k}{\l+m}\left[kyj_{\vk-1}(ky)\right]^2 \right)dy, \qq\l\in\s_{\rm ac}(H_0),\]
and $k=k(\l)=\sqrt{\l^2-m^2}$ is defined in (\ref{q-m}). Here $j_{\vk}$ is spherical Bessel function given  by (\ref{jk}).

Moreover,
\[\lb{Ome0}\Omega(\l)=\Omega_0 +{\mathcal O}\left(\frac{\ln |\l|}{|\l|}\right)\qq\mbox{as}\,\,\l\rightarrow\infty,\qq \mbox{where}\,\,\Omega_0=\int_0^\infty v(x)dx.\]

The function $\phi_{\rm sc}$ (defined in (\ref{Omega})) is continuous on
$\R\sm [-m,m]$ and satisfies
\[
\lb{asf}
\phi_{\rm sc}(\l)=\O_0+{\mathcal O}\left(\frac{\ln |\l|}{|\l|}\right)\qq\mbox{as}\qq \l \rightarrow \pm\infty.
\]

 If in addition $v'\in L^1(\R_+),$ then the functions $\gf^+(\l),  D(\l)$   satisfy for  $\l\in\C_+$
\[
\lb{a=D} \gf^+(z)=k_0(z) D(z)\exp\left(i\Omega_0+\frac{1}{\pi}\int\frac{\Omega(t)-\Omega_0}{t-z}dt\right). 
\]
\end{theorem}
{\bf Remark.}  The condition $v\in L^\infty([0,a])$ for some $a>0$ is needed in (\ref{Ome0}- \ref{a=D}).

\subsection{Resonances}
In order to consider  resonances we need a stronger hypothesis on the function $v$.

\no {\bf Condition A}. {\it Real-valued function $v\in L^2(\R_+)$
and $\supp v\ss [0,\g], \g>0$,
 where  $\g=\sup\supp v.$ }

Later we will even suppose that, in addition, $v'\in L^1(\R_+).$

 We denote $\sqrt{z}$ the principal branch of the square root that
is positive for $z>0$ and with the cut along the negative real axis.
We denote $\C_\pm=\{\l\in\C; \pm\Im\l >0\}.$ 

We introduce the
quasi-momentum $ k(\l)$ by
\[\lb{q-m}
k(\l)=\sqrt{\l^2-m^2},\qqq \l\in \L=\C\sm [-m,m].
\]
The function $ k(\l)$ is a conformal mapping from $\L$ onto
$\cK=\C\sm [im,-im]$ and satisfies
\[
k(\l)=\l-{m^2\/2\l}+{{\mathcal O}(1)\/\l^2}\qq \mbox{as} \qq |\l|\to \iy.
\]
The function $k(\l)$ maps the horizontal cut  $(-m,m)$ on the
vertical cut $[im,-im]$. Moreover,
$$
k(\R_\pm \sm (-m,m))=\R_\pm,\qqq k(i\R_\pm)=i\R_\pm \sm (-im,im).
$$
The Riemann surface for $k(\l)$ is obtained  by joining the upper
and  lower rims of two copies  $\C\sm \s_{\rm ac}(H_0)$ cut along
the $\s_{\rm ac}(H_0)$ in the usual (crosswise) way. Instead of this
two-sheeted Riemann surface it is more convenient to work on the cut
plane $\L$ and half-planes $\L_\pm$ given by
$$
\L=\C\sm [-m,m],\qqq \L_\pm=\C_\pm\cup g_\pm.
$$
Here we denote $g_+\subset\L_+,$  and $g_{-}\subset\L_-,$ the upper
respectively and lower rim of the cut $(-m,m)$ in $\C\setminus
[-m,m].$ Here the upper half-plane $\L_+=\C_+\cup g_+$ corresponds
to the physical sheet and the lower half-plane $\L_+=\C_-\cup g_-$
corresponds to the non-physical sheet.

{\bf Below we consider all functions and the resolvent in $\C_+$ and
will obtain their analytic continuation throught the continuous spectrum $\s_{\rm ac}(H_0)$ into the cut domain $\L$.}

Note that, equivalently, we could consider the Jost function, the
resolvent etc in $\L_-$ (the physical sheet) and obtain their
analytic continuation into the whole cut domain $\L$.

{\bf By abuse of notation, we will think of all functions $f$ as functions of both $\l$ and $k,$ and will regard notations as $f(\cdot,\l),$  $f(\cdot,k)$ and similar as indistinguishable.}

It is well known that for each  $h\in C_0(\R_+,\C^2)$ the function  $\Phi(\l)=((H-\l)^{-1} h,h)$ has meromorphic continuation from  $\C_+$ into $\C\setminus \s_{\rm ac}(H_0).$ We denote $g^{+}\subset\L_1^+,$  $g^{-}\subset\L_1^-,$ the upper respectively lower rim of the gap $(-m,m)$ in $\C\setminus [-m,m].$

\begin{definition} \lb{DefRes} Let $\Phi(\l)=((H-\l)^{-1} h,h),$ $\l\in\C\setminus \s_{\rm ac}(H_0)$ for some $h\in C_0(\R_+,\C^2),$ $h\neq 0.$\\
\no 1) If $\Phi(\l)$ has pole at some $\l_0\in g^+$  we call $\l_0$ an {\bf eigenvalue}.\\
\no 2) If $\Phi(\l)$ has pole at some $\l_0\in\L_1^-$  we call $\l_0$ a {\bf resonance}.\\
\no 3) A point $\l_0 =m$ or $\l_0=- m$ is called {\bf virtual state} if the function $z\rightarrow \Phi(\l_0 +z^2)$ has a pole at $0.$\\
\no 4) A point $\l_0\in\L$  is called a {\bf state} if it is either an eigenvalue, a resonance or a virtual state. Its multiplicity is the multiplicity of the corresponding pole. We denote $\s_{\rm st}(H)$ the set of all states. If $\l_0\in \s_{\rm st}(H)\cap g^{-},$ then we call $\l_0$ an {\bf anti-bound state}.
\end{definition}

We will show that the set of resonances coincides with the set of zeros in $\L_1^-$ of the Jost function $\gf^+(\l)$ defined in (\ref{Jost_function_def}) or, equivalently, of the modified Fredholm determinant $D(\l).$ 
Multiplicity of a resonance is the multiplicity of the corresponding zero.

Recall that $\vk$ is the spin-orbit coupling parameter defined in (\ref{defvk}).

\begin{proposition}\lb{freevirt}
For $\vk\neq 0$ operator  $H_0$ does not have virtual states.
\end{proposition}
{\bf Remark.}  In \cite{IantchenkoKorotyaev2013a} it was shown that in regular case (which corresponds to  $\vk=0$) the point $\l=-m$ is the virtual state of $H_0.$

We show that the following results valid in the regular case as in \cite{IantchenkoKorotyaev2013a} also hold in framework of the present paper.

\begin{theorem}\lb{Th-bound-antibound} Let $V$ satisfy condition A. Then the states of $H$ satisfy:\\
1) The number of eigenvalues is finite.\\
2) Let $\l^{(1)}\in g^{+}\subset\L_1^+$ be eigenvalue of $H$ and $\l^{(2)}\in g^{-}\subset\L_1^-$ be the same number but on the "non-physical sheet". Then $\l^{(2)}$ is not an anti-bound state.\\
3) Let $\l_1,\l_2\in g^{+}, $ $\l_1 <\l_2,$ be  eigenvalues of $H$ and assume that there are no other eigenvalues on the interval $\o^{(1)}=(\l_1,\l_2)\ss g^{+}.$ Let $\o^{(2)}\ss g^{-}$ be the same interval but on the "non-physical sheet". Then there exists an odd number $\geq 1$ of anti-bound states (counted with multiplicities) on $\o^{(2)}.$
\end{theorem}

{\bf Remark.} Resonances for operators with gaps in the spectrum were studied
in \cite{Korotyaev2011}, \cite{KorotyaevSchmidt2012}. In these papers the investigation
of the resonances  on the cut plane was transformed into the theory
of the entire function theory. For the Dirac operators on the half-plane with
$\vk=0$ it was done in \cite{IantchenkoKorotyaev2013a}. In our paper we use
similar arguments.

\

An entire  function $f(z)$ is said to be
of exponential type  if there is a constant $A$ such that
$|f(z)|\leq\const e^{A|z|}$ everywhere. The infimum of the set of
$A$ for which inequality holds is called the type of $f(z)$ (see
\cite{Koosis1988}). Section 2 in in \cite{IantchenkoKorotyaev2013a}  contains more details on the
exponential type functions. If $f$ is analytic and satisfies the
above inequality only in $\C_+$ or $\C_-,$ we will say that $f$ is
of exponential type in $\C_\pm$ with the type defined appropriately. Applying  a version of Froese Lemma \ref{l-Fr} we get in Theorem  \ref{Th-ass2terms}   the exact exponential type of the Jost functions. 

This result allows to determine the asymptotics of the counting function. We denote the
number of zeros of a function $f$ having modulus  $\leq r$ by $\cN
(r,f)$, each zero being counted according to its multiplicity.

\begin{theorem}
\lb{Prop_counting_zeros}  Let the potential $v$ satisfy Condition A and $v'\in L^1(\R_+).$ Then $D(\cdot)$ has  an analytic extension from $\C_+$ into
the whole cut plane $\C\setminus [-m,m].$  The set of zeros of $D$  satisfy:
\[
\lb{counting}
 \cN(r,D)={2r\g\/ \pi }(1+o(1))\qqq as \qqq r\to\iy.
\]
For each $\d >0$ the number of zeros of $D$ with negative imaginary part with modulus $\leq r$
lying outside both of the two sectors $|\arg z |<\d ,$ $|\arg z -\pi
|<\d$ is $o(r)$ for large $r$.

\end{theorem}
{\bf Remark.} 1) Zworski  obtained in \cite{Zworski1987} similar results for the Schr{\"o}dinger operator with compactly supported potentials on the real line.

2) Our proof follows from Proposition \ref{Prop_counting_zeros_F} and the Levinson Theorem.
\vspace{3mm}

\subsection{Trace formulas}
In the massless case $m = 0$ we have $ k = \l,$  $k_0 = −i,$ $\sigma_{\rm ac}(H_0) = \R,$ and the Riemann surface consists of two disjoint sheets $\C$ (see \cite{IantchenkoKorotyaev2013} for the regular case). Therefore we can consider all functions and the resolvent in
the upper-half plane $\C_+$ and  obtain their analytic
continuation into the whole complex plane $\C$. Then the Jost functions are analytic on $\C$ and allows Hadamard factorization (\ref{Had-fact}).

Denote by $\{\l_n\}_{1}^{\iy}$
the sequence of its zeros of the Jost function $\gf^+$ in $\C_-$ (multiplicities counted by repetition), i.e. complex resonances,
so arranged that $0<|\l_1|\leq |\l_2|\leq |\l_2|\leq \dots.$

We prove the following theorem.

\begin{theorem}
\lb{T4} Let the mass $m=0$ and let the potential $v$ satisfy Condition A. Let
$f\in\mS$ where $\mS$ is the Schwartz class of rapidly decreasing functions. Let $\l_n$ denote either a resonance if
$\l_n\neq 0$ or an eigenvalue if $\l_n=0.$ Let $\f_{\rm sc}(\l)$ be the scattering phase.
Then
\[\lb{TrF1}
\Tr (f(H)-f(H_0))=-\frac{1}{\pi}\int_\R f(\l)\f_{\rm sc}'(\l)d\l, \]
\[\lb{sc_phase}
\f_{\rm sc}'(\l)=\gamma+\sum_{|\l_n|\neq 0}\frac{\Im\l_n}{|\l-\l_n|^2}, \qq \l\in\R.\]
\[\lb{TrF3}
 \Tr (R(\l)-R_0(\l))=-i\g-\lim_{r\to +\infty} \sum_{|\l_n|\le
r}\frac{1}{\l-\l_n},\qq \Im\l\neq 0,
\]
where the series converge uniformly in every bounded subset on the
plane by condition (\ref{sumcond}).
\end{theorem}

{\bf Remark.}  Such identities were obtained for Schr\"odinger operators 
on the half-line \cite{Korotyaev2004} and  were  extended to massless Dirac operators (regular case) in \cite{IantchenkoKorotyaev2013}. 
 In our paper we use
similar arguments.

{\bf The plan of paper is as follows.} In Section \ref{s-Free Dirac system} we collect all needed facts related to the  unperturbed radial Dirac operator $H_0.$ The  proof of Proposition
\ref{freevirt}  is given there. Moreover, we remind the associated spectral representation and study the Hilbert-Schmidt norms of the cut-off free resolvent, Theorem \ref{L-HS} and traces,
Lemma \ref{res2}.

 In Section \ref{s-AsJost} we define and study the Jost functions as well as we prove Theorem \ref{Th-ass2terms} using a version of Froese Lemma \ref{l-Fr}. In order to achieve this goal we will need to get uniform estimates on the Jost function.  

In Section \ref{s-MFD} we give the properties of the modified Fredholm determinant and prove  the main result of the paper Theorem \ref{T1}. Moreover, we give a useful expression for the trace  of the difference of the resolvents, Proposition \ref{P-RR}.

In Section \ref{s-F} we introduce and study an analytic function $\gF$ which is used in order to prove Theorems \ref{Th-bound-antibound} and \ref{Prop_counting_zeros}.

In Section \ref{s-massless} we study the massless case and prove the trace formulas stated in Theorem \ref{T4}.

We moved to  the appendix, Section \ref{S-proofL-2}, the (quit technical) proof of Lemma \ref{L-2}.

\section{ Free Dirac system.}\lb{s-Free Dirac system}
\setcounter{equation}{0}

\subsection{Preliminaries}
For  the free radial Dirac operator $H_0$ we consider the corresponding
free radial Dirac system
\[
\lb{Dirac_system_0}
 \ca      f_1'+\frac{\vk}{x}f_1-(m+\l)f_2=0 \\
           f_2'-\frac{\vk}{x}f_2-(m-\l)f_1=0, \ac \ \ \ \l \in \C,\qqq 
           f= \ma f_1\\ f_2\am,
\]
where $f_1, f_2$ are  the functions of $x\in\R_+ $. System (\ref{Dirac_system_0}) can be written equivalently as
\[
\lb{Dsx}     f'=\ma -{\vk\/x} & m+\l\\ m-\l & {\vk\/x}\am 
f\qq\Leftrightarrow \qqq
           \ma f_2 \\   f_1 \am '=\ma {\vk\/x} & m-\l\\ m+\l & -{\vk\/x}\am 
           \ma f_2 \\    f_1 \am.
\]

 Recall that
 $k=k(\l)=\sqrt{\l^2-m^2}$ and $k_0(\l)=\frac{\l +m}{i k(\l)}.$

We consider the  fundamental solutions $\vp, \vt$ of (\ref{Dirac_system_0})  satisfying
 $\det (\vt,\vp)=1$ and
\[
 \lb{vtat0}
 \vp(x,\l)={x^{\vk}\/(2\vk-1)!!}\ma 0\\ 1\am (1+o(1)),\qqq
 \vt(x,\l)={(2\vk-1)!!\/x^{\vk}}\ma 1\\ 0\am (1+o(1)),
 \] 
as $ x\to 0, \l\in \C$,  where $(2\vk-1)!!=1\cdot 3\cdot 5\cdot\ldots\cdot (2\vk-1),$ if $(2\vk-1)\geq 1,$ and $(2\vk-1)!!=1,$ if $(2\vk-1)\leq 0.$
Thus $\vp, \vt$ are  given by
\[
\lb{vp0def}
\vp(x,\l)=\ma\vp_1(x,\l) \\
\vp_2(x,\l) \am=k^{-\vk}\ma \frac{k(\l)}{\l-m}zj_\vk(z)\\
zj_{\vk-1}(z)  \am=k^{-\vk}\ma ik_0(\l)\si\\
\co  \am ,
\]
\[
\lb{vt0def}
\vt(x,\l)=\ma\vt_1(x,\l) \\
\vt_2(x,\l) \am=k^{\vk}\ma z\eta_\vk(z)\\
\frac{k(\l)}{\l+m}z\eta_{\vk-1}(z) \am ,\qq z=kx,
\]

Here $j_\vk(z)$ is the Spherical Bessel function of the first kind,
$$
zj_\vk(z)=\sqrt{\pi z\/2} J_{\vk+{1\/2}}(z),
$$
and $J_\nu$ is the Bessel function (see \cite{Erdelyietal1953II}, p.4 formula 2).
Moreover, note the following useful formula
\[\lb{jk}
zj_\vk(z)=\frac{\pi^\frac12}{2}z\sum_{\ell=0}^\infty\frac{(-1)^\ell}
{\ell!\G(\ell+\vk+\frac32)}\left({z\/2}\right)^{2\ell+\vk},
\]
which implies that $\si$ is odd if $\vk$ is even, and  $\si$ is even if $\vk$ is odd, the property which will be used later in this paper.

Now, we introduce a basis of {\bf Jost solutions}  $\psi^\pm$ for (\ref{Dirac_system_0})
  \[\lb{JostAsy}
 \begin{aligned}& H_0\psi^\pm=\l\psi^\pm,\qq\psi^\pm(x,\l)=(\mp i k)^\vk e^{\pm i k(\l) x}\ma\pm k_0(\l)\\ 1\am+o(1),\,\,\mbox{as}\,\,x\rightarrow\infty,\\ & k_0(\l)=\frac{\l +m}{i k(\l)},\,\, \l\in\sigma_{\rm ac}(H_0).
\end{aligned}
  \]
Using (\ref{JostAsy}), the Wronskian for the pair $\psi^+,\psi^-$ is then given by \[\lb{Wronskian_free} 
\det \left(\psi^+(\cdot,\l),\psi^-(\cdot,\l)\right)=-2i(\l+m)k^{2\vk-1}.
\] 
The Jost solutions are represented using the Spherical Bessel functions as follows
\[\lb{freeJostsol}
\psi^\pm(x,\l)=\mp i k^\vk \ma {\l+m\/k} zh_\vk^\pm(z) \\
                        zh_{\vk-1}^\pm (z)  \am=\mp i k^\vk                       \ma ik_0(\l) zh_\vk^\pm(z) \\
                        zh_{\vk-1}^\pm (z)  \am,\qq z=kx,
\] 
with  
$$
zh_\vk^\pm(z)=z(\eta_\vk(z)\pm ij_\vk(z)),
$$ 
and where
$$-z\eta_\vk(z) =  z y_\vk(z)=(-1)^\vk \sqrt{\pi z\/2} J_{-\vk-{1\/2}}(z).$$
Here $y_\vk(z)$ is the Spherical Bessel function of the second kind. We have also the following relations:
$$
z y_\vk(z)=\sqrt{\pi z\/2} Y_{\vk+{1\/2}}(z)=(-1)^{\vk +1} zj_{-\vk-1}(z),
$$
 $$
 zh^+_0(z)=e^{iz}, \ zh^+_{-1}(z)=ie^{iz}, \ zj_0(z)=\sin z,\qq zj_{-1}(z)=\cos z.
 $$
 We have
\[
\lb{Hankel}\begin{aligned}
 zh_\vk^+(z)=i\sqrt{\pi z\/2}H_{\vk+{1\/2}}^{(1)}(z)=i\sqrt{\pi z\/2}\left(J_{\vk+{1\/2}}(z)+iY_{\vk+{1\/2}}(z)\right)
=iz\left(j_\vk(z)+iy_\vk(z)\right),\\
 zh_\vk^-(z)=-i\sqrt{\pi z\/2}H_{\vk+{1\/2}}^{(2)}(z)=-i\sqrt{\pi z\/2}\left(J_{\vk+{1\/2}}(z)-iY_{\vk+{1\/2}}(z)\right)
=-iz\left(j_\vk(z)-i y_\vk(z)\right),
\end{aligned}
\]
and the Hankel functions $H_\nu^{(j)}(z),$ $j=1,2,$ are defined in \cite {Erdelyietal1953II}, p.4, formulas 5, 6.
Asymptotics  from \cite{Erdelyietal1953II}, 7.13.1, page 85, imply the asymptotics for $|z|\rightarrow\infty$
$$
\begin{aligned}&zh_\vk^+(z)=i\left(\frac{\pi z}{2}\right)^{1/2}H_{\vk+1/2}^{(1)}(z)\sim e^{i(z-\frac{\pi}{2}\vk)},\\
 &zh_\vk^-(z)=-i\left(\frac{\pi z}{2}\right)^{1/2}H_{\vk+1/2}^{(2)}(z)\sim e^{-i(z-\frac{\pi}{2}\vk)},
\end{aligned}$$
and therefore asymptotics (\ref{JostAsy}) for Jost functions $\psi^\pm.$

We collect some useful formulas in the two lemmas below.  The proof follows from \cite{Erdelyietal1953II} (see also \cite{Barthelemy1967} and \cite{Blancarteetal1995}).

\begin{lemma}\lb{l-bounds on fundsol}  Let $\vk\in\Z_+$.  Then uniformly in $z=xk\in \C\setminus \{0\}$ the following estimates hold true:
 \[\lb{2.13}
  |z j_\vk(z)|\leq C e^{|\Im z|}\left(\frac{|z|}{1+|z|}\right)^{\vk+1},
\]
\[\lb{2.14}
|z \eta_\vk(z)|\leq C e^{|\Im z|}\left(\frac{1+|z|}{|z|}\right)^{\vk},
\]
\[\lb{bound_vp}
|\vp(x,\l)|\leq   Ce^{|\Im z|}\left(\frac{x}{1+|z|}\right)^{\vk}\ma \frac{|\l+m|x}{1+|z|}\\ 1  \am ,
\]
\[\lb{bound_psi}
|\psi^\pm (x,\l)|\leq   C|k| e^{\mp \Im z }\left(\frac{1+|z|}{x}\right)^{\vk-1}\ma \frac{1+|z|}{|\l-m|x}\\ 1  \am ,
\]

\[\lb{bound_vpTpsi} |\vp^T(x,\l)\psi^+(x,\l)|\leq C e^{(|\eta|-\eta)x}\frac{|z|}{1+|z|}\left(\left|\frac{\l+m}{\l-m}\right|+1\right)\leq C' e^{(|\eta|-\eta)x}
\]
\[\lb{bound_vpTpsi2} |\vp^T(x,\l)\psi^+(x,\l)-1|\leq C e^{(|\eta|-\eta)x}\left(\frac{1}{1+|z|}+\frac{1}{|k|}\right).
\]
\end{lemma}

\begin{lemma}
\lb{l-bounds on fundsol} 
Let $\F=kx-\vk\frac{\pi}{2}, z=xk$. For each $\vk\in\Z_+$ the following asymptotics hold true:\\
for $|z|\rightarrow\infty,$
\[\lb{vp_exp} 
\vp(x,\l)= k^{-\vk}\rt\{\ma  ik_0(\l)\sin \F \\
                  \cos\F \am+\ma ik_0(\l)\frac{\vk(\vk+1)}{2z}\cos\F \\
\frac{-\vk(\vk-1)}{2z}\sin \F \am
 +\cO \left({e^{|\Im z|}\/|z|^2}\right)\rt\},
\]
\[\lb{psi_exp} 
\psi^+ (x,\l)=k^\vk e^{i\F}\rt\{\ma k_0(\l) \\ 1 \am+\ma
ik_0(\l) \frac{\vk(\vk+1)}{2z} \\
i \frac{\vk(\vk-1)}{2z} \am
+{\cO(1)\/|z|^2}\rt\};
\]
for $|z|\to 0,$
\[
 \lb{besselat0}%Barthelemy:
zj_\vk(z)=\frac{(z)^{\vk+1}}{(2\vk+1)!!}+{\mathcal O}(z^{\vk+3}),\]
 \[
 \lb{bessel2at0}
 z\eta_\vk(z)=z^{-\vk}(2\vk -1)!!+{\mathcal O}(z^{-\vk+2}), 
 \]
\[\lb{ash} 
zh_\vk^\pm(z)=\frac{(2\vk-1)!!}{z^\vk}+{\mathcal O}(z^{-\vk+2}).
\]

\end{lemma}

Now, we have the following representation of the Jost solution 
\[\lb{Jost_Weyl} 
\psi^+(x,\l)=k^{2\vk} \left[ k_0(\l)k^{-2\vk}\vt(x,\l)  + \vp(x,\l)\right],
\] 
as
$$
\begin{aligned}
&i\psi^+(x,\l)= k^\vk\ma \frac{\l+m}{k} zh_\vk^+(z) \\
zh_{\vk-1}^+ (z)  \am
=k^\vk\ma  \frac{\l+m}{k} z\eta_\vk(z) \\
z\eta_{\vk-1}(z)  \am+ik^\vk\ma  \frac{\l+m}{k} zj_\vk(z) \\
                        zj_{\vk-1}(z)  \am=\\ 
&= ik^{2\vk}\left[\frac{1}{i}k^{-\vk}\left(
                      \begin{array}{c}
                        \frac{\l+m}{k} z\eta_\vk(z) \\
                        z\eta_{\vk-1}(z)  \\
                      \end{array}
                    \right)+
                    k^{-\vk}\left(
                      \begin{array}{c}
                        \frac{\l+m}{k} zj_\vk(z) \\
                        zj_{\vk-1}(z)  \\
                      \end{array}
                    \right) \right]=\\ \\
&=i k^{2\vk} \left[\frac{\l+m}{ik}k^{-2\vk}k^{\vk}\left(
                      \begin{array}{c}
                         z\eta_\vk(z) \\
                        \frac{k}{\l+m}z\eta_{\vk-1}(z)  \\
                      \end{array}
                    \right)+
                    k^{-\vk}\left(
                      \begin{array}{c}
                        \frac{\l+m}{k} zj_\vk(z) \\
                        zj_{\vk-1}(z)  \\
                      \end{array}
                    \right) \right]=\\ \\
&=  ik^{2\vk} \left[ k_0(\l)k^{-2\vk}\vt(x,\l)  + \vp(x,\l)\right].
\end{aligned}
$$

This yealds  the free  radial Titchmarch-Weyl function
 \[\lb{Weyl}
 m_\vk= \frac{(k(\l))^{2\vk}}{k_0(\l)}=i\frac{(k(\l))^{2\vk+1}}{\l+m}.
 \]

We define the {\bf Jost function} $\gf^{0,+}(\l)$  for the unperturbed radial Dirac operator by
\[\lb{Free_Jost_function_def}
\gf^{0,+}(\l)=\det \left(\psi^+(\cdot,\l),  \vp(\cdot,\l)\right)=\lim_{x\to 0}\frac{x^\vk}{(2\vk-1)!!}\psi_1^+(x,\l).\]

Now, using (\ref{freeJostsol}) and (\ref{ash}) we get in the
 leading order
\[\lb{psiat0}\psi^+(x,\l)\sim -ik^\vk \left(
                                 \begin{array}{c}
                                   \frac{\l+m}{k}\frac{(2\vk-1)!!}{z^\vk} \\
                                   \frac{(2\vk-3)!!}{z^{\vk-1}} \\
                                 \end{array}
                               \right),\qq\mbox{as}\,\, |k|x\rightarrow 0\]
and
$$
\gf^{0,+}(\l)=\lim_{x\rightarrow 0}\frac{x^\vk}{(2\vk-1)!!}(-i)k^\vk\frac{\l+m}{k}\frac{(2\vk-1)!!}{z^\vk}
=(-i)\frac{\l+m}{k}=k_0(\l).
$$

The conjugate Jost function $\psi^-(\l)$ is then given by  
$$
\psi^-(\l)=-\det (\vp(\cdot,\l),\psi^-(\cdot,\l))=-k_0(\l),\qq \l\in\s(H_0),
$$ 
which yelds
\[\lb{relcon}\begin{aligned}
\vp(x,\l)&=\frac{i}{2(\l+m) k^{2\vk-1}}\left(-\psi^-( \l)\psi^+(x,\l)+\psi^+( \l)\psi^-(x,\l)\right)\\&=\frac{1}{2 (k(\l))^{2\vk}}\left(\psi^+(x,\l)+\psi^-(x,\l)\right).
\end{aligned}\]

{\bf Resolvent.}
Now, the integral kernel of the free resolvent $R_0(\l):=(H_0 -\l)^{-1}$ is given by
$$ R_0(x,y,\l)=\left\{\begin{array}{lr}
                 \frac{1}{\det(\psi^+,\vp)}\psi^+(x,\l)(\vp(y,\l))^T & \mbox{if}\,\, y<x, \\
                 \frac{1}{\det(\psi^+,\vp)}\vp(x,\l)(\psi^+(y,\l))^T & \mbox{if}\,\, x<y.
               \end{array}\right.
$$
Using (\ref{psiat0}) and (\ref{vtat0}) we get $\det\left(\psi^+(\cdot,\l),\vp(\cdot, \l,\vk)\right)=\frac{\l+m}{ik(\l)}=k_0(\l).$

Then,
$$
R_0(x,y,\l)=\ma
        \frac{k(\l)}{\l-m}zh^+_\vk(z)\, \z j_\vk(\z) & zh^+_\vk(kx)\,\z j_{\vk-1}(\z) \\ zh^+_{\vk-1}(z)\,\z j_\vk(\z) & \frac{k(\l)}{\l+m}zh^+_{\vk-1}(z)\,\z j_{\vk-1}(\z) \am,\qq z=kx,\,\,\z=ky,\qq\mbox{ if}\,\,y<x.
    $$
\vspace{3mm}

{\bf Proof of Proposition \ref{freevirt}. }

We have
\[\lb{prop_k_m}
\begin{aligned}
&k(\l)=i\sqrt{2m}\sqrt{\epsilon}\left( 1-{\mathcal O}(\epsilon)\right),\qq k_0(\l)=-\frac{\sqrt{2m}}{\sqrt{\epsilon}}\left( 1+{\mathcal O}(\epsilon)\right),\qq\epsilon=m-\l\rightarrow 0+,\\
&k(\l)=i\sqrt{2m}\sqrt{\epsilon}\left( 1-{\mathcal O}(\epsilon)\right),\qq k_0(\l)=-\frac{\sqrt{\epsilon}}{\sqrt{2m}}\left( 1+{\mathcal O}(\epsilon)\right),\qq\epsilon=m+\l\rightarrow 0+.
\end{aligned}
\]

Formulas (\ref{besselat0})  (\ref{ash}) for $\vk=1,2,\ldots$ imply
$$\vp(x,\l)\sim x^{\vk}\ma (\l+m)\frac{x}{ (2\vk+1)!!}\\ \\ \frac{1}{ (2\vk-1)!!}\am,\qq\psi^+(x,\l)\sim \frac{-i}{x^\vk} \left(
                                 \begin{array}{c}
                                   \frac{\l+m}{k}(2\vk-1)!! \\ \\
                                   (2\vk-3)!!(kx) \\
                                 \end{array}
                               \right),\qq kx\rightarrow 0,
$$ 
we get for $y<x$
 $$
 R_0(x,y,\l)=\frac{1}{k_0(\l)}\psi^+(x,\l)(\vp(y,\l))^T\sim \frac{k(\l)}{\l+m}\left(
 \begin{array}{cc}                                                                            \frac{(\l+m)^2}{k}\frac{x(2\vk-1)!!}{(2\vk+1)!!} & \frac{\l+m}{k(\l)} \\                                                                             (\l+m)\frac{x}{ (2\vk+1)!!}  (2\vk-3)!!(kx)& \frac{kx}{2\vk-1} \\                                                                          \end{array}                                                                        \right).
 $$

Therefore, $\{\pm m\}$ are not virtual states for $\vk\geq 1.$\hfill\BBox

 {\bf Remark.} The situation for $\vk=0$ is different (see \cite{IantchenkoKorotyaev2013a}). For $\vk=0$ there is virtual state  $\{-m\}$ . This does not contradict our proof for $\vk\geq1$ as if $\vk=0$ we can not use the asymptotics at zero which shows that $ \frac{kx}{2\vk-1}\rightarrow 0,$ as $k\rightarrow 0.$  Instead, for $\vk=0$ we should write   $kxh^+_{-1}(kx)=ie^{ikx}\rightarrow i,$ $k\rightarrow 0,$ which would imply that $\{-m\}$ is virtual state for $\vk=0.$

\subsection{ Spectral representation} In this section we follow the classical ideas of spectral representation for Dirac operators \cite{LevitanSargsyan1988} as presented in \cite{IantchenkoKorotyaev2013a}  in the regular case.

Let, as before,$$\vp(x,\l)=\ma\vp_1(x,\l) \\
\vp_2(x,\l) \am=k^{-\vk}\ma \frac{k(\l)}{\l-m}\si\\
\co  \am $$ be the regular at $x=0$ solution. Then there  exists a non-decreasing function $\rho(s),$ $s\in\R,$ such that
 for any vector-function $f\in L^2(\R_+)$  there exists function $\cF\in L^2(\R,d\rho)$ such that
 \[\lb{spectral representation}
 \cF(s)=\int_0^\infty f^{\rm T}(x)\vp(x,s)dx,\qq f(x)=\int_{-\infty}^\infty \cF(s)\vp(x,s)d\rho(s).
 \]
and
\[\lb{Parseval} \int_0^\infty(f_1^2(x)+f_2^2(x))dx=\int_{-\infty}^\infty \cF^2(s)d\rho(s).\]
Function $\rho$ is the spectral function. It satisfies the finiteness condition $\int_{-\infty}^\infty (1+s^2)^{-\vk-1}d\rho(s) <\infty.$
Here, $\cF$ is the generalized Fourier transform of the vector-function $f$ with respect to the solutions of the Dirac equation (\ref{Dirac_system_0}) with the Dirichlet boundary condition. We denote the generalized Fourier transform by $\Phi$ and write $\cF(s)=(\Phi f)(s).$ Formula (\ref{Parseval}) is the Parseval's identity and it shows that $\Phi$ is an isometry of $L^2(\R_+)\times(L^2(\R_+)$ onto $L^2(\R,d\rho).$

As the discreet spectrum of $H_0$ is empty, then $\rho(s)=0$ for $s\in (-m,m).$ For $\l\in\sigma_{\rm ac}(H_0)$
the function $\rho$ can be easily derived from the Weyl function $m_\vk(\l)$ obtained in (\ref{Weyl}).

For $s\in \sigma_{\rm ac}(H_0)=(-\infty,-m]\cup[m,+\infty)$ we get, using (\ref{Weyl}),  $$d\rho(s)=\rho'(s)ds=\frac{1}{\pi}\Im m_\vk(s+i0)ds=\frac{1}{\pi}\frac{(k(s))^{2\vk+1}}{s+m}ds.$$

Using that $\rho'(s)=\frac{1}{\pi}\frac{(k(s))^{2\vk+1}}{s+m}$ is positive  for $\l\in(-\infty,-m]\cup[m,\infty)$ and by introducing  the functions
$\hat{\vp}(x,s)=\vp(x,s)\sqrt{\rho'(s)},$ $\hat{\cF}(s)=\cF(s)\sqrt{\rho'(s)},$ we get
\[\lb{spectral representation2}
 \hat{\cF}(s)=\int_0^\infty f^{\rm T}(x)\hat{\vp}(x,s)dx,\qq f(x)=\int_{-\infty}^{-m} \hat{\cF}(s)\hat{\vp}(x,s)ds+\int_m^\infty \hat{\cF}(s)\hat{\vp}(x,s)ds.
 \]
and
\[\lb{Parseval2} \int_0^\infty(f_1^2(x)+f_2^2(x))dx=\int_{-\infty}^{-m} \hat{\cF}^2(s)ds+\int_m^\infty \hat{\cF}^2(s)ds.\]

The modified generalized Fourier transform  $\hat{\Phi}:$ $\hat{\cF}(s)=(\hat{\Phi} f)(s)$  is an isometry of $\cH=(L^2(\R_+)^2$ onto
$$\hat{\cH}=\hat{\Phi}(\cH)=L^2((-\infty, -m],\rho'(s)ds)\oplus L^2([m,+\infty),\rho'(s)ds).$$
Moreover, as for any $f\in\cH,$ $g\in\hat{\cH},$ we have
$$\langle\hat{\Phi}f,g\rangle_{\hat{\cH}}=\langle f, \hat{\Phi}^{-1}g\rangle_{\cH}$$ and
$\hat{\Phi}^{-1}$ is the formal adjoint $(\hat{\Phi})^*$ of $\hat{\Phi}.$ Here $\langle\cdot,\cdot\rangle_{H}$ denotes the scalar product in the Hilbert space $H.$

Let $$\cE(x,s)=\left(
           \begin{array}{cc}
             \hat{\vp}_1(x,s) & 0 \\
             0 & \hat{\vp}_2(x,s) \\
           \end{array}
         \right)$$ and $\sigma:=\s_{\rm ac}(H_0)=(-\infty, -m]\cup[m,\infty).$ Then
it follows from (\ref{Parseval2})

\[\lb{Parseval3} \int_0^\infty\left|\int_\s \hat{\cF}(s)\cE(y,s)ds\right|^2dy=\int_\s |\hat{\cF}(s)|^2ds \]

As $$\hat{\vp}(x,s)=\sqrt{\frac{1}{\pi}\frac{(k(s))^{2\vk+1}}{s+m}}(k(s))^{-\vk}\ma \frac{k(s)}{s-m}\si\\
\co  \am =\frac{1}{\sqrt{\pi}}\left(
                      \begin{array}{c}
                       \sqrt{\frac{s+m}{k(s)}}\si \\
                       \sqrt{\frac{k(s)}{s+m}}\co \\
                      \end{array}
                    \right),$$ then using
the bound (\ref{2.13})
$$\left|\si\right|\leq Ce^{|\Im k|x}\left( \frac{|k|x}{1+|k|x}\right)^{\vk+1}$$ applied for $\l\in\sigma\subset\R,$ we get
 $$\begin{aligned}\pi |\cE(x,s)|^2=&\left|\frac{s+m}{k(s)}\right|(\si)^2 +\left|\frac{k(s)}{s+m}\right|(\co)^2 \\
\leq& \max \left( \left|\frac{s+m}{k(s)}\right|,\left|\frac{k(s)}{s+m}\right|\right)Ce^{2|\Im k(s)|x}\left( \frac{|k(s)|x}{1+|k(s)|x}\right)^{2\vk}\leq C\cK_1(s),\end{aligned}$$
where $$\cK_1(s)=\max \left( \left|\frac{s+m}{k(s)}\right|,\left|\frac{k(s)}{s+m}\right|\right).$$
\subsection{ Hilbert-Schmidt norms}

We define the sets \[\lb{set}
\cZ_\epsilon^\pm=\{\l\in\C\setminus [-m,m];\,\,\pm\Im\l\geq 0,\,\,|\l\pm m|\geq \epsilon\},\qq\cZ_\epsilon=\cZ_\epsilon^+\cup\cZ_\epsilon^-,\qq\epsilon>0.
\]

We denote by $\|.\|_{\cB_k},$ the trace ($k=1$) and the Hilbert-Schmidt ($k=2$) operator norms.

For a Banach space $\cX$ let $AC(\cX)$ denote the set of all
$\cX$-valued analytic functions on $\C_+,$ continuous in $\overline{\C}_+\setminus\{\pm m\}.$

\begin{theorem}
\lb{L-HS} Let $\chi,\widetilde{\chi}\in L^2(\R_+;\C^2)$ and $\l\in \C\setminus\R.$
Put \[\lb{CL} C_\lambda=\left[
\frac{4\pi {\mathcal O}(1)}{|\Im\l|}\left|\Re\frac{\l}{\sqrt{\l^2-m^2}}\right|+{\mathcal O}\left(\max\left\{\frac{1}{|\l|^2},\frac{1}{|\l\pm m|^2}\right\}\right)\right].
\]
Then it follows:\\
i) Operators $\chi R_0(\l), R_0(\l)\chi, \chi R_0(\l) \wt\chi$ are the $\cB_2$-valued
functions satisfying the following properties:
\[
\begin{aligned}
\lb{B2-1} \|\chi R_0(\l)\|_{\cB_2}^2=\| R_0(\l)\chi\|_{\cB_2}^2 \leq C_\lambda   \|\chi\|_2^2,
\end{aligned}
\]
\[
\lb{B2-2} \begin{aligned}&\|\chi R_0(\l)\wt\chi \|_{\cB_2}\le \frac{c}{\epsilon}\|\chi\|_2\|\wt\chi\|_2\qq\mbox{for}\,\,\l\in\cZ_\epsilon,\\ &\|\chi R_0(\l)\wt\chi \|_{\cB_2}\rightarrow 0\,\,\mbox{as}\,\,|\Im\l|\rightarrow\infty.\end{aligned}
\]
Moreover, for each $\l\in \C_+,$ the operator-function $\chi R_0\wt\chi\in AC(\cB_2)$.

ii) For each $\l\in \C\setminus\R,$ operator $ \chi R_0'(\l) \wt\chi=\chi R_0^2(\l) \wt\chi\in AC( \cB_2)$ is the $\cB_2$-valued
function satisfying
\[
\lb{B2-3}
\|\chi R_0'(\l)\wt\chi\|_{\cB_2}\le \frac{c}{\epsilon^2}\|\chi\|_2\|\wt\chi\|_2\qq \mbox{for}\,\,\l\in \cZ_\epsilon,\qq \|\chi R_0'(\l)\wt\chi \|_{\cB_2}\rightarrow 0\,\,\mbox{as}\,\,|\Im\l|\rightarrow\infty.\]
\end{theorem}

The proof of Theorem \ref{L-HS} is identical to that in the regular case and is given in full detail in \cite{IantchenkoKorotyaev2013a}.  It is based on spectral representation of the resolvent via the
  generalized Fourier  transform $\Phi.$ Here we will repeat only some arguments which will be also used later (Lemma \ref{res2}).

{\bf Proof.}
Let $\sigma=\s_{\rm ac}(H_0)=(-\infty, -m]\cup[m,\infty)$ and $\Im\l\neq 0.$
 The transformed free resolvent acting in the  $s-$space  $L^2(\R,d\rho(s))$  is the operator of multiplication
 by $\frac{1}{s-\l}$ and we have
 \[\lb{L-HS-spec}\begin{aligned}
 R_0(\l)f(x)=&\int_{-\infty}^\infty \frac{1}{s-\l}\left(\int_0^\infty f^{\rm T}(t)\vp(t,s)dt\right)\vp(x,s) d\rho(s)\\
 =&\int_\s \frac{1}{s-\l}\left(\int_0^\infty f^{\rm T}(t)\hat{\vp}(t,s)dt\right)\hat{\vp}(x,s) ds=\\
 =&\int_\s R_0(x,s,\l)\int_0^\infty \cE(y,s)f(y)  dy ds
 ,\end{aligned}
\] where
\[\lb{kernel_in_fourier}\begin{aligned}
&
 R_0(x,s,\l)=\frac{1}{s-\l}\left(
                                                                 \begin{array}{cc}
                                                                   \hat{\vp}_1(x,s)   & 0 \\
                                                                   0 & \hat{\vp}_2(x,s)   \\
                                                                 \end{array}
                                                               \right)= \frac{1}{s-\l}\cE(x,s),\\ &\cE(x,s)=\left(
           \begin{array}{cc}
             \hat{\vp}_1(x,s) & 0 \\
             0 & \hat{\vp}_2(x,s) \\
           \end{array}
         \right).\end{aligned}\]
Let $\chi\in L^2(\R_+;\C^2)$ and we can suppose that $\chi$ is diagonal matrix. Then using  (\ref{Parseval3})
$$\int_0^\infty\left|\int_\s \hat{\cF}(s)\cE(y,s)ds\right|^2dy=\int_\s |\hat{\cF}(s)|^2ds$$ we get
\[\lb{121}
\begin{aligned}
\|\chi R_0(\l)\|_{\cB_2}^2&=\int_0^\infty\!\!\int_0^\infty\left|\int_\s \chi (x) R_0(x,s,\l)\cE(y,s)ds\right|^2dy dx=\int_0^\infty\int_\s |\chi (x) R_0(x,s,\l)|^2dsdx\\
&=\int_0^\infty\int_\s \left|\chi (x)M(x,s)\frac{1}{s-\l}\right|^2dsdx\leq \frac{C}{\pi}\int_0^\infty|\chi(x)|^2dx\int_\s \frac{\cK_1(s)}{|s-\l|^2}ds,
\end{aligned}
\] where
$$ \cK_1(s)=\left( \left|\frac{s+m}{k(s)}\right|,\left|\frac{k(s)}{s+m}\right|\right).$$
The rest of the proof  is identical to the regular case $\vk=0$ and can be found in \cite{IantchenkoKorotyaev2013a}.\hfill\BBox
\vspace{3mm}

In order to proof the trace formula we will need the following lemma which follows directly from the spectral representation of the resolvent (\ref{L-HS-spec}).
\begin{lemma}\lb{res2}  Let $V=vI_2\in L^2(\R_+;\C^2).$ For any $\l\in\C\setminus\R$  $$\Tr (VR_0'(\l))= \Tr (VR_0^2(\l))=\int_0^\infty\!\!\int_\s\frac{v(x)}{(s-\l)^2}\left( \hat{\vp}^2_1(x,s)+\hat{\vp}^2_2(x,s)\right) ds dx=\frac{1}{\pi}\int_\sigma\frac{\O(s)}{(s-\l)^2}ds,$$ where $\Omega(s),$ $s\in\s=\s_{\rm ac}(H_0),$ is given by 
\[\lb{Oid}
\Omega(\l)=\int_0^\infty v(y)\left( \frac{k(\l)}{\l-m}\left[kyj_\vk(ky)\right]^2 +\frac{k(\l)}{\l+m}\left[kyj_{\vk-1}(ky)\right]^2 \right)dy,\qq \l\in\s_{\rm ac}(H_0).
\]
\end{lemma}

\section{Asymptotics of the Jost solutions}\lb{s-AsJost}
\setcounter{equation}{0}

The main goal of this section is to get uniform estimates on the Jost solutions needed in order to get exact exponential type of the Jost functions,  Theorem \ref{Th-ass2terms}. 

\subsection{Preliminaries}
For  the radial Dirac operator $H$ we consider the corresponding
radial Dirac system:
\[\lb{Dirac_system}
 \ca  f_1'+\frac{\vk}{x}f_1-(m-v(x)+\l)f_2=0 \\
           f_2'-\frac{\vk}{x}f_2-(m+v(x)-\l)f_1=0\ac,\ \ \ \l \in \C,\qqq
           f= \ma f_1\\ f_2\am.
\]
 Note that for any two solutions $f,g $ of (\ref{Dirac_system})  the Wronskian $\det(f,g)=f_1g_2-f_2g_1$ is independent of $x.$
The regular case $\vk=0$ was studied in \cite{IantchenkoKorotyaev2013a}.
For $\vk\neq 0$ the problem (\ref{Dirac_system})
 is singular at $x=0.$

 We consider the regular solution $\phi(x,\l)$ of (\ref{Dirac_system})
 satisfying
 \[
 \lb{phi0}
 \f(x,\l)=\frac{x^{\vk}}{(2\vk-1)!!}\ma 0\\ 1\am \qqq \mbox{as} \qq x\to 0.
 \]
We introduce the {\bf Jost solutions}  $f^\pm$ for (\ref{Dirac_system}) by the  conditions 
\[
\begin{aligned}
  \lb{JostAsy1}
&  Hf^\pm=\l f^\pm,\qq f^\pm(x,\l)= (\mp i k)^\vk e^{\pm i k(\l) x}\ma\pm k_0(\l)\\ 1\am+o(1)\qq
  \mbox{as}\,\,x\rightarrow\infty,
  \\ & k_0(\l)=\frac{\l +m}{i k(\l)},\qq \l\in\sigma_{\rm ac}(H_0),
\end{aligned}
\]
 where $k=k(\l)=\sqrt{\l^2-m^2}$ was defined in (\ref{q-m}). Note that $f^-(x,\l)=\overline{f^+(x,\l)}$ for $\l\in\sigma_{\rm ac}(H).$
Recall that the Jost solutions for the unperturbed system ($v=0,$ associated with free radial Dirac operator (\ref{free_radial})) are defined by the same conditions (\ref{JostAsy}) and are denoted by $\psi^\pm(x,\l).$

Using the regular solution $\phi$ with asymptotics (\ref{phi0}) we define the {\bf Jost function} by
\[\lb{Jost_function_def}
\gf^+(\l)=\det \left(f^+(\cdot,\l),  \phi(\cdot,\l)\right)=\lim_{x\rightarrow 0}\frac{x^\vk}{(2\vk-1)!!}f_1^+(x,\l).
\]
We denote $\gf^{0,+}(\l)$  the Jost function for
 the unperturbed Dirac system ($v=0$). Recall that   $\gf^{0,+}(\l)=k_0(\l).$

{\bf Remark.} Our definition of the Jost solutions by asymptotics (\ref{JostAsy1}) implies that
$$f^+(x,\l)=k_0(\l)(\theta(x,\l) +m(\l)\phi(x,\l)),$$ where $\theta,\phi$ are the fundamental solutions of (\ref{Dirac_system}) satisfying
 $\det(\theta,\phi)=1$ and (\ref{phi0}); $m=m(\l)$ is the Titchmarsch-Weyl function, which is equal to $$m_\vk=\frac{(k(\l))^{2\vk}}{k_0(\l)}=i\frac{(k(\l))^{2\vk+1}}{\l+m}$$
 in the unperturbed case $v\equiv 0.$ The choice of normalization in  (\ref{JostAsy1}) implies also that the unperturbed Jost function $\gf^{0,+}(\l)$ is independent of $\vk$ and is the same as in the regular case $\vk=0$ discussed in \cite{IantchenkoKorotyaev2013a}.

Using  asymptotics (\ref{JostAsy1}) we get
that the Wronskian of the pair $f^+,f^-$ is given by
 \[
 \lb{wrid}
\det (f^+,f^-)=\det (\psi^+,\psi^-)=2k_0k^{2\vk},
\] 
where we used (\ref{Wronskian_free}). 

\

The main result of this section is the following theorem 
\begin{theorem}\lb{Th-ass2terms} Assume that the potential $v$ satisfies Condition A and   $v'\in
L^1(\R).$ Then the Jost function $\gf^\pm(\cdot)$ has exponential type $2\gamma$ in $\C_\mp.$
\end{theorem}

\

 In order to prove Theorem \ref{Th-ass2terms} we need to study analytic properties of the Jost functions. We start with deriving the integral equation for the Jost solution.

Let $\vt(\cdot)=\vt(\cdot,\vk,\l),$  $\vp(\cdot)=\vp(\cdot,\vk,\l)$ be fundamental solutions to $H_0$ with $\det(\vt,\vp)=1,$ satisfying  $$ \lim_{x\rightarrow 0}x^{-\vk}\vp(x,\vk,\l)=\frac{1}{(2\vk-1)!!}\ma 0\\ 1\am,\qq \lim_{x\rightarrow 0}x^{\vk}\vt(x,\vk,\l)=(2\vk-1)!!\ma 1\\ 0\am.$$

 Then $f^+$ satisfies the integral equation
$$
\begin{aligned}
 f^+(x,\l)=\psi^+(x,\l)+\int_x^\infty
 G(x,y,\l)V(y)f^+(y,\l)dy,\\
 G(x,y,\l)=-\left(\vp(x)\vt^{\rm T}(y)-\vt(x)\vp^{\rm T}(y)\right).
\end{aligned}
$$
 Using asymptotics (\ref{Jost_function_def}),  (\ref{vtat0}) we get that the Jost functions satisfies
$$
\gf^+(\l)=k_0(\l)+\int_0^\infty v(y)
 \vp^{\rm T}(y)f^+(y,\vk,\l)dy=k_0+\int_0^\infty v(y)
 \left(\vp_1f^+_1+\vp_2f^+_2\right)dy.
 $$

Put
$$U(x)=\left(
      \begin{array}{cc}
        \frac{1}{\psi_1^+(x)} & 0 \\
        0 & \frac{1}{\psi_2^+(x)}  \\
      \end{array}
    \right)=\left(
      \begin{array}{cc}
        \frac{1}{k^\vk k_0 kxh_\vk^+(kx)} & 0 \\
        0 & \frac{1}{-ik^\vk  kxh_{\vk-1}^+(kx)}  \\
      \end{array}
    \right),\qq \chi=Uf.$$
Then $\chi$ satisfies
$$
\c(x)=\c^0
+\int_x^\infty U(x)G(x,y)U^{-1}(y)V(y)\chi(y)dy, \qq \c^0=\ma 1\\ 1\am,
$$
as $U^{-1}$ and $V$ commute. Thus we have the power series
\[\label{series}
\c (x,\l )=\sum _{n\geq 0}\c^n(x,\l ),\ \ \ \
\c^{n+1}(x,\l )=\int_x^\infty U(x)G(x,y)U^{-1}(y)V(y)\c^n(y,\l )dy.
\]

We formulate the following standard result without a proof. The first part of Lemma \ref{chi_estimates} was shown in \cite{Blancarteetal1995}, the proof of the second part is straightforward.
This Lemma is generalization for the singular potential $q=\vk/x$ of Lemma 4.1 in
\cite{IantchenkoKorotyaev2013a}.

\begin{lemma}\lb{chi_estimates} Let $\e :=\Im\,k(\l)$ and $M$ the  matrix valued function
$$ \cM(x,y,\l)=\left(
                                                                                                                                         \begin{array}{cc}
                                                                                                                                           |k_0| C\left(\frac{1+|k|x}{|k|x}\right)^{2\vk}\left(\frac{1+|k|y}{|k|y}\right)^{2\vk}  & \frac{C}{|k_0|} \left(\frac{1+|k|x}{|k|x}\right)^{2\vk} \left(\frac{1+|k|y}{|k|y}\right)^{2(\vk-1)} \\
                                                                                                                                           C|k_0| \left(\frac{1+|k|x}{|k|x}\right)^{2(\vk-1)}\left(\frac{1+|k|y}{|k|y}\right)^{2\vk} & \frac{C}{|k_0|}\left(\frac{1+|k|x}{|k|x}\right)^{2(\vk-1)}\left(\frac{1+|k|y}{|k|y}\right)^{2(\vk-1)} \\
                                                                                                                                         \end{array}
                                                                                                                                       \right)$$
 and for each $\l\in\cZ_\epsilon^+$ let  $M(\l,\delta)=\sup_{y>x>\delta}\|\cM(x,y,\l)\|,$

1) Suppose $v\in L^1(\R_+)$ and  $x\geq\delta>0,$  $\epsilon >0.$ Then  the function $\chi(x,\cdot)$ is analytic in $\cZ_\epsilon^+,$ and for 
$\l\in \cZ_\epsilon^+,$ the functions $\chi^n,$ $\chi$ satisfy the following estimates:
\[\lb{xn+}
  \|\c^n(x,\l )\|\leq \frac{1}{n!}\left(M(\l,\delta)\int_x^\infty |v(t)|dt\right)^{n},\ \ \forall \ n\geq 1,\]
\[
\lb{x+}
  \|\chi(x,\l)\|\leq e^{ M(\l,\delta) \int_x^\infty |v(t)| dt}.
\]

2) If $v$ satisfies Condition A, then for each $x\in \R_+,$ $\epsilon >0$  the function $\chi(x,\cdot)$ is analytic in $\C\setminus\{\pm m\}.$
For each $x\in [0,\gamma],$ $\epsilon >0$
and $\l\in \cZ_\epsilon,$ the vector functions $\chi^n,$ $\chi$ satisfy the following estimates :

\[\lb{xn}
  \|\c^n(x,\l )\|\leq e^{(\gamma-x)(|\e|-\e)}\frac{1}{n!}\left(M(\l,\delta)\int_x^\gamma |v(t)|dt\right)^{n},\ \ \forall \ n\geq 1,\]
\[
\lb{x}
  \|\chi(x,\l)\|\leq e^{(\gamma-x)(|\e|-\e)}e^{ M(\l,\delta) \int_x^\gamma |v(t)|dt}.
\]

\end{lemma}

From this Lemma it follows
\begin{corollary}\lb{Cor-jostsolutions} Let $v\in L^1(\R_+,\C^2)$ and  $x\geq\delta >0,$  $\epsilon >0.$
1)  Then  the function $f^+(x,\cdot)$ is analytic in $\cZ_\epsilon^+.$

2) If, in addition $v$ satisfies Condition A, then  the function $f^+(x,\cdot)$ is analytic in $\C\setminus\{\pm m\}.$

\end{corollary}

We recall the following results (see Theorem 3.1 in \cite{Blancarteetal1995}).

\begin{theorem}\lb{th-BFW1} Let $v\in L^1(\R_+).$ Then for $\l\in\C_+,$ as $|k|\rightarrow\infty$ and $|k|x\rightarrow\infty,$  the following facts hold true:\\
\[\lb{i)}f^+(x,\l)= e^{i\int_x^\gamma v(t)dt}\psi^+(x,\l)\left[1+(1+x) {\mathcal O}\left(\frac{1}{|k|x}\right)\right],\]
 \[\lb{ii)}f^+(x,\l)= e^{i\int_x^\gamma v(t)dt}\psi^+(x,\l)\left[1+\left(1+\frac{1}{x}\right)\frac{g(k)}{|k|x}+ {\mathcal O}\left(\frac{1}{|k|^2}\right)\right], \]
where $g(k)=o(1),$ and if also $v'\in L^1,$  then $g(k)={\mathcal O}(|k|^{-1})$ and
 \[\lb{iii)} f^+(x,\l)= e^{i\int_x^\gamma v(t)dt}\psi^+(x,\l)\left[1+\left(1+x\right)\frac{C}{|k|^2x^2}+ {\mathcal O}\left(\frac{1}{|k|^2}\right)\right].\]
\end{theorem}
\begin{theorem}[Theorem 3.2 in \cite{Blancarteetal1995} ]\lb{th-BFW2} Suppose $v\in L^1(\R_+)$ and for some $a>0,$ $1\leq q <\infty,$ $v\in L^q(0,a).$ Then as $|\l|\rightarrow\infty,$ $\l\in\C_+,$
$$\gf^+(\l)=-i e^{i\int_0^\infty v(t)dt}+o\left(\frac{1}{|k|^{(q-1)/q}}\right).$$
If $v\in L^1\cap L^\infty(0,a),$ then
$$\gf^+(\l)=-i e^{i\int_0^\infty v(t)dt}+{\mathcal O}\left(\frac{\ln |k|}{|k|}\right).$$
\end{theorem}

\vspace{5mm}

Let $\phi(x,\l)$ be regular solution satisfying (\ref{phi0}). Then for $\l\in\sigma_{\rm ac}(H_0)$ we get (compare with the free case (\ref{relcon}))
$$\phi(x,\l)=\frac{i}{2(\l+m) k^{2\vk-1}}\left(-f^-(\l)f^+(x,\l)+f^+( \l)f^-(x,\l)\right).$$ In the limit $ x\rightarrow \infty,$  $\phi(x,\l)$ behaves asymptotically as 
$$ \frac{i}{2(\l+m) k^{2\vk-1}}\left(-f^-( \l)\psi^+(x,\l)+f^+( \l)\psi^-(x,\l)\right),$$ and 
if, in addition, $v$ satisfies Condition A, two functions coincide for $x\geq\gamma.$

We write
$$\gf^+(\l)=|\gf^+(\l)|e^{i\delta_\vk (\l)},\qq \gf^-(\l)=\overline{\gf^+(\l)},\qq\l\in\sigma_{\rm ac}(H_0).$$
%Then as our Jost solution is related to
%Klaus by $f^+=k_0\widetilde{f}^+$ and $\arg k_0(\l)=-\pi/2$ for $\l >m,$ then we get $\delta_\vk=\widetilde{\delta}_\vk-\pi/2.$
Using asymptotics (\ref{JostAsy1}) we get in the main order as  $ x\rightarrow \infty$
$$\phi(x,\l)\sim\frac{|\gf^+(\l)|}{(\l+m) k^{\vk-1}}\left(
                                                         \begin{array}{c}
                                                           \frac{\l+m}{k(\l)}\cos \left(kx-\vk\frac{\pi}{2}-\delta_\vk (\l)\right) \\
                                                           \sin \left(kx-\vk\frac{\pi}{2}-\delta_\vk (\l)\right) \\
                                                         \end{array}
                                                       \right).$$

\subsection{Uniform estimates on the Jost solutions}\lb{s-uniformJost}
In order to get uniform estimates on the Jost function as $|\l|\rightarrow\infty$ we need to transform the Dirac system (\ref{Dirac_system}) to more convenient form. We follow \cite{IantchenkoKorotyaev2013a}. To start with, the free (non-radial) Dirac equation
$$\left(
    \begin{array}{c}
      f_1 \\
      f_2 \\
    \end{array}
  \right)'=\left(
            \begin{array}{cc}
              0 & m+\l \\
              m-\l & 0 \\
            \end{array}
          \right)
  \left(
            \begin{array}{c}
              f_1 \\
              f_2 \\
            \end{array}
          \right)$$ is transformed to the diagonal form
$$\left(
    \begin{array}{c}
      \widetilde{f}_1 \\
      \widetilde{f}_2 \\
    \end{array}
  \right)'=\left(
            \begin{array}{cc}
              ik(\l) & 0 \\
              0 & -ik(\l) \\
            \end{array}
          \right)
  \left(
            \begin{array}{c}
              \widetilde{f}_1 \\
              \widetilde{f}_2 \\
            \end{array}
          \right),\,\, \widetilde{f}=Uf, $$
          where  $$U=\frac12\left(
                              \begin{array}{cc}
                                \frac{1}{\l+m} & \frac{1}{ik(\l)}  \\
                                \frac{1}{\l+m}  & -\frac{1}{ik(\l)} \\
                              \end{array}
                            \right),\qq U^{-1}=\left(
                                                 \begin{array}{cc}
                                                   \l+m & \l+m \\
                                                   ik(\l) & -ik(\l) \\
                                                 \end{array}
                                               \right)=ik(\l)\left(
                                                 \begin{array}{cc}
                                                  k_0 & k_0\\
                                                   1 & -1 \\
                                                 \end{array}
                                               \right).$$

                                               If $f=\psi^+$ is the unperturbed radial Jost function, we get
                                               $\widetilde{f}(x,\l)=U\psi^+(x,\l),$ where $\psi^+$ is defined in (\ref{freeJostsol})
$$ 
\psi^+(x,\l)=- i k^\vk\ma \frac{\l+m}{k} zh_\vk^+(z) \\
  zh_{\vk-1}^+ (z)  \am,\qq zh_\vk^\pm(z)=z(\eta_\vk(z)\pm ij_\vk(z)),\qq z=kx .
  $$ 
  Thus
$$
\widetilde{f}(x,\l)=-{i k^{\vk -1}\/2} \left(
                      \begin{array}{c}
                         z(h_\vk^+(z)-  ih_{\vk-1}^+ (z)) \\
                        z(h_\vk^+(z) +  ih_{\vk-1}^+ (z))  \\
                      \end{array}
                    \right)\sim -i k^{\vk -1} \left(
                      \begin{array}{c}
                         e^{i(kx-\frac{\pi}{2}\vk)} \\
                       0 \\
                      \end{array}
\right),\qq |z|\rightarrow\infty,
$$
where by using Formula 3, page 78 in \cite{Erdelyietal1953II}, we get
$$\begin{aligned}
&z(h_\vk^+(z)-  ih_{\vk-1}^+ (z))=(-i)^\vk e^{iz}\left(\sum_{j=0}^{\vk-1}i^j\frac{(\vk +j-1)!}{j!(\vk-j)!}\underline{2\vk} (2z)^{-j}+i^\vk\frac{(2\vk)!}{\vk !}(2z)^{-\vk}\right)\\
&z(h_\vk^+(z)+  ih_{\vk-1}^+ (z))=(-i)^\vk e^{iz}\left(\sum_{j=0}^{\vk-1}i^j\frac{(\vk +j-1)!}{j!(\vk-j)!}\underline{2j} (2z)^{-j}+i^\vk\frac{(2\vk)!}{\vk !}(2z)^{-\vk}\right).
\end{aligned}
$$
We write 
$$
\widetilde{f}(x,\l)=-i k^{\vk -1} e^{i(kx-\frac{\pi}{2}\vk)}  \left(
\begin{array}{c} 1+A \\   B \\ \end{array}   \right),
$$ 
where for any $\delta >0$  the functions (where $z=kx$)
\[\lb{AB}\begin{aligned}
&A=A(z,\vk)= \frac12\left(\sum_{j=1}^{\vk-1}i^j\frac{(\vk +j-1)!}{j!(\vk-j)!}\underline{2\vk} (2z)^{-j}+i^\vk\frac{(2\vk)!}{\vk !}(2z)^{-\vk}\right)=\sum_{j=1}^\vk a(\vk,j)z^{-j},\\
& 
B=B(z,\vk)=\frac12\left(\sum_{j=1}^{\vk-1}i^j\frac{(\vk +j-1)!}{j!(\vk-j)!}\underline{2j} (2z)^{-j}+i^\vk\frac{(2\vk)!}{\vk !}(2z)^{-\vk}\right)=\sum_{j=1}^\vk b(\vk,j)z^{-j}
\end{aligned}\]
are uniformly bounded in any bounded sub-domain
of $\{(k,x)\in \cZ_\epsilon^+\times\R_+;\qq |kx|>\delta\}$ and for any $\vk\in\N,$ $A,B\rightarrow 0$ as $|kx|\rightarrow\infty.$

Now, put
\[\lb{Mtl} 
\mM=\ma  0 & \ol {M}\\  M & 0  \am ,\qq M(t,\l)=e^{-i2\int_0^t v(s)ds}\left(\frac{\vk}{t}-\frac{im}{k}v(t)\right),\]
\[\lb{Ntl}  
\mN=\ma  N & 0\\ 0 & \ol{N}  \am,\qq  N(t,\l)=\frac{i(\l-k)}{k}v(t)={\mathcal O}(k^{-2}).
\]

\begin{lemma}\lb{L-1} Suppose $v$ satisfy Condition A.

 Let $f^+$ be the Jost solution and let the vector-function $Y=Y(x,\l)$ be defined via 
 $$
 f^+(x,\l)=k^\vk e^{i\left(v_0-\frac{\pi}{2}\vk\right)}
 \ma  k_0(\l) & k_0(\l) \\  1 & -1  \am 
 \ma    e^{-i\int_0^t v(x)dx} & 0 \\    0 &  e^{i\int_0^t v(x)dx} \am   Y.
 $$
where $v_0=\int_0^\g v(x)dx$.
Then
$Y=Y(x)$ satisfies the differential equation
\[\lb{odeY}Y'(x)=\left(ik\s_3-(\mN+\mM)\right)Y,\qq Y(x,\l)_{|x\geq\g}=e^{ikx}\left(
                        \begin{array}{cc}
                          1 & 0 \\
                          0 & e^{-i2v_0} \\
                        \end{array}
                      \right)\left(
                                                                                             \begin{array}{c}
                                                                                               1+A \\
                                                                                               B \\
                                                                                             \end{array}
                                                                                           \right),\]
which is equivalent to the integral equation
\[\lb{inteqY}Y(x,\l)=e^{ikx}\ma 1 & 0 \\ 0 & e^{-i2v_0} \am 
\ma   1+A \\ B  \am         +\int_x^\gamma e^{ik\sigma_3(x-t)}\left(\mN(t,\l)+\mM(t,\l )\right)Y(t,\l)dt.\]
\end{lemma}
{\bf Proof.} Firstly,  similar to \cite{IantchenkoKorotyaev2013a} and originally \cite{HKS88}, \cite{HKS89}, by a chain of transformations of the Dirac equation (we omit the details here), we introduce a new vector-function $X$ related to the Jost solution $f^+$ via
 \[\lb{fX}f^+=ik\left(
                                                            \begin{array}{cc}
                                                              k_0 & k_0 \\
                                                              1 & -1 \\
                                                            \end{array}
                                                          \right)e^{i\sigma_3\left(kx-i\int_0^x v(t)dt\right)}X.
                                                   \]
Then $X$ satisfies the differential equation
\[\lb{diffeqX}\begin{aligned}&X'=-\wt{W}X,\qq X_{|x\geq\g}=X^0:= -ik^{\vk-1}e^{i(kx-\frac{\pi}{2}\vk)}\left(
\begin{array}{cc}
e^{-ikx+i\int_0^\g v(t)dt} & 0 \\
0 & e^{ikx-i\int_0^\g v(t)dt} \\
\end{array}
\right)\left(\begin{array}{c}
    1+A \\
        B \end{array}\right),\\& \widetilde{W}(t)=e^{-ikt\s_3}\left(\mN+\mM\right)e^{ikt\s_3}.
\end{aligned}
\]
Then $\widetilde{X}=i k^{-\vk +1} e^{-i\left(\int_0^\gamma v(t)dt-\frac{\pi}{2}\vk\right)} X(t)$ satisfies 
$$
\wt{X}'=-\wt{W}\wt{X},\qq\wt{X}(x)=\left(
                      \begin{array}{c}
                        1+A \\
                        0 \\
                      \end{array}
                    \right) + e^{i2[kx-\int_0^\gamma v(t)dt]}\left(
                    \begin{array}{c}
                     0 \\                                                                                             B \\                                                                                             \end{array}                                                                                           \right) +\int_x^\infty\wt{W}(t)\wt{X}(t)dt.
$$

 We write $$\left(
                     \begin{array}{cc}
                       N & \overline{M}\\
                       M & \overline{N}  \\
                     \end{array}
                   \right)=\left(
                     \begin{array}{cc}
                       N & 0\\
                       0 & \overline{N}  \\
                     \end{array}
                   \right)+\left(
                     \begin{array}{cc}
                       0 & \overline{M}\\
                       M & 0  \\
                     \end{array}
                   \right)=\mN+\mM.$$

                    We introduce new vector-valued function $Y$ by $$Y=\left(
                     \begin{array}{cc}
                       e^{ikt} & 0 \\
                       0 & e^{-ikt}  \\
                     \end{array}
                   \right)\widetilde{X}=e^{ikt\sigma_3}\widetilde{X},\qq \widetilde{X}= \left(
                     \begin{array}{cc}
                       e^{-ikt} & 0 \\
                       0 & e^{ikt}  \\
                     \end{array}
                   \right)Y= e^{-ikt\sigma_3}Y. $$

Now, the function $$Y(x,\l)=ik^{-\vk+1}e^{ikx\s_3}e^{-i\left(\int_0^\g v(s)ds-\frac{\pi}{2}\vk\right)}X(x,\l)$$ satisfies (\ref{odeY}) and (\ref{inteqY}).\hfill\BBox

\begin{lemma}\lb{L-2} 
Let $v$ satisfy Condition A and in addition $ v'\in L^1(\R_+;\C).$
We denote  $$\mW=\mW(t,\l)=2ik\mN+\mM'(t)-\mM\mN-|M|^2,\qq \mA= \mA(x,\l)=I-\frac{1}{2ik}\sigma_3\mM.$$
The function $\mW(t,\l)$ has the following asymptotics as $|k|\rightarrow\infty:$
$$\begin{aligned}
&\mW(t,\l)=\mW_0(t,\l)+{\mathcal O}(k^{-1}),\qq\mW_0(t,\l)=\ma -\frac{\vk^2}{t^2}& \ol w(t,\l) \\ w(t,\l) & -\frac{\vk^2}{t^2}\am, \\ & w(t,\l)=-e^{-i2\int_0^t v(s)ds}\left[\frac{\vk}{t^2}+i2v(t)\frac{\vk}{t}+\frac{i(\l-m)\vk}{kt}v(t)\right]
\end{aligned}
$$

Then for $|kx|\geq \sup_{x\in\R_+}|xM(x)|$ the matrix $\mA(x,\l)=I-\frac{1}{2ik}\sigma_3\mM(x,\l)$ has bounded inverse $\mB(x,\l)=I+{\mathcal O}\left( (|k|x)^{-1}\right)$ and the solution  $Y=Y(x,\l)$ of equations (\ref{odeY}), (\ref{inteqY})
 satisfies
 $$\begin{aligned}&Y=Y^0+(2ik)^{-1}\mB KY,\qq Y^0=Y^0(x,\l)=e^{ikx}\mB(x,k)\left(
                        \begin{array}{cc}
                          1 & 0 \\
                          0 & e^{-i2\int_0^\gamma v(t)dt} \\
                        \end{array}
                      \right)\left(
                                                                                             \begin{array}{c}
                                                                                               1+A \\
                                                                                               B \\
                                                                                             \end{array}
                                                                                           \right), \\& KY=\int_x^\gamma
e^{i\sigma_3k(x-t)}\mW(t,\l)Y(t,\l)dt,\end{aligned},$$ where functions $A,B$ are defined in (\ref{AB}),    and $Y=Y(x,\l)$ is given by the expansion in powers of $(2ik)^{-1}$
$$Y=Y^0+\sum_{n\geq 1}Y^n,\qq  Y^n=\frac{1}{(2ik)^n}(\mB K)^nY^0,$$ where
 $$|Y^n(x,k)|\leq \frac{2C_{\epsilon,\delta}}{n!|k|^n}e^{|\Im k|(2\gamma-x)}\left(\int_x^\gamma |\mW(s)|ds\right)^n,  $$
  and $$C_{\epsilon,\delta}=\sup_{(k,x)\in \Omega_{\epsilon,\delta}}\{|1+A(kx)|,|B(kx)|\},\qq \Omega_{\epsilon,\delta}=\{(k,x)\in \cZ_\epsilon^+\times\R_+;\qq \min\{ |k|x, x\}>\delta\}.$$

Moreover,
$$Y^0=e^{ikx}\left[\left(
                     \begin{array}{c}
                       1 \\
                       0 \\
                     \end{array}
                   \right)+(1+x){\mathcal O}\left(\frac{1}{|k|x}\right)
\right],\qq Y^n=e^{|\Im k|(2\gamma-x)}{\mathcal O}\left(\frac{1+x}{|k|x}\right)^n.$$

\end{lemma}

{\bf Remark.} The "Moreover" statement follows from $\mW(t,\l)=t^{-2}{\mathcal O}(1).$

The proof of Lemma \ref{L-2} is given in Section \ref{S-proofL-2}. Here, we will apply this Lemma in order to prove Theorem \ref{Th-ass2terms}.

We will need the following Lemma by Froese (see  \cite{Froese1997}, Lemma 4.1). Even though the original lemma was stated for  $U\in L^\infty$ the argument also works for $U\in L^2$ and we omit the proof.
\begin{lemma}[Froese]\lb{l-Fr}  Suppose $U\in L^2(\R)$ has compact support
contained in $[0,1],$ but in no smaller interval. Suppose $g(x,\l)$
is analytic for $\l$ in the lower half plane, and for real $\l$ we
have $g(x,\l)\in L^2([0,1]\,dx,\R\, d \l).$ Then $\int e^{i\l x}
U(1+g(x,\l))\,d x$ has exponential type at least $1$ for $\l$ in the
lower half plane.\end{lemma}

{\bf Proof of Theorem \ref{Th-ass2terms}.} We will use the following relations
$$f^+(x,\l)=k^\vk e^{i\left(\int_0^\gamma v(t)dt-\frac{\pi}{2}\vk\right)}\left(
                                                                   \begin{array}{c}
                                                                     k_0e^{-i\int_0^x v(t)dt}Y_1+k_0e^{i\int_0^x v(t)dt}Y_2  \\
                                                                      e^{-i\int_0^x v(t)dt}Y_1-e^{i\int_0^x v(t)dt}Y_2 \\
                                                                   \end{array}
                                                                 \right)
                                                      ,$$
$$\gf^+(\l)=k_0+\int_0^\infty v(y)
 \left(\vp_1f^+_1+\vp_2f^+_2\right)dy,\qq \vp=k^{-\vk}\ma ik_0\,\si\\
\co  \am,$$ which implies $$\begin{aligned}
\gf^+(\l)=&k_0+e^{i\left(\int_0^\g v(t)dt-\frac{\pi}{2}\vk\right)}k^\vk\int_0^\infty v(y)\left[e^{-i\int_0^y v(t)dt}Y_1(y)\left(k_0\vp_1(y)+\vp_2(y)\right)\right.\\
&+\left. e^{i\int_0^y v(t)dt}Y_2(y)\left(k_0\vp_1(y)-\vp_2(y)\right)\right]dy.
\end{aligned} $$

Let $\widetilde{X}(y)=e^{-ik y\sigma_3}Y(y).$ Then, for $\l\in \C_-,$  from the properties of $Y$ as in the proof of  Lemma \ref{L-2} it follows that $\widetilde{X}_1=1+g(y,k),$ $g(y,k)=(1+y){\mathcal O}\left(\frac{1}{|k|y}\right),$ and $\widetilde{X}_2=(1+y){\mathcal O}\left(\frac{1}{|k|y}\right).$

Put $\phi_\pm:=e^{-iky}\left(k_0\vp_1(y)\pm \vp_2(y))\right).$ For $\l\in \C_-,$ we have (see proof of (\ref{2.13}))  $$|\phi_\pm|=\left|e^{-iky}\left(k_0\vp_1(y)\pm \vp_2(y))\right)\right|\leq C \left(\frac{|k|x}{1+|k|x}\right)^{\vk+1}.$$

We write  $$\begin{aligned}
\gf^+(\l)=&k_0+e^{i\left(\int_0^\g v(t)dt-\frac{\pi}{2}\vk\right)}k^\vk\int_0^\infty v(y)e^{2iky}\left[e^{-i\int_0^y v(t)dt}\widetilde{X}_1(y)\phi_+(y)+ e^{i\int_0^y v(t)dt}\widetilde{X}_2(y)\phi_-\right]dy\\
=&k_0+e^{i\left(\int_0^\g v(t)dt-\frac{\pi}{2}\vk\right)}k^\vk\int_0^\infty v(y)\phi_+(y)e^{2iky}e^{-i\int_0^y v(t)dt}\left(1+g(y,k)\right)dy\\
& +e^{i\left(\int_0^\g v(t)dt-\frac{\pi}{2}\vk\right)}k^\vk\int_0^\infty v(y)e^{2iky} e^{i\int_0^y v(t)dt}\widetilde{X}_2(y)\phi_-dy
\end{aligned}$$

Let $K(\l)=\int_0^\g e^{2iky} U(y)\left(1+g(y,k)\right)dy,$ $U(y)=v(y)\phi_+(y)e^{-i\int_0^y v(t)dt}.$
Now, it is enough to apply a version \ref{l-Fr} of Lemma of Froese to $K(\l-i),$ $\l\in\C_-,$  where we shift the argument of function $K$ in order to avoid the singularities at $\l=\pm m,$ and using that $\sup_{t\in [0,\gamma]} | g(t,k(\tau -i))|={\mathcal O}(\tau^{-1})$ as $\tau\rightarrow\pm\infty.$ Thus the function $f_1(0,k)$ has exponential type $2\gamma$ in the half plane $\C_-,$ and Theorem \ref{Th-ass2terms} is proved.\hfill\BBox

\section{Modified Fredholm determinant}\lb{s-MFD}
\setcounter{equation}{0}

The main goal of this section is to prove the main result of the paper --  Theorem \ref{T1}.

In order to prove Theorem \ref{T1} to study the properties of the modified Fredholm resolven.

Recall definitions in (\ref{Y0}). $R_0(\l)=(H_0 -\l)^{-1},$ $R(\l)=(H -\l)^{-1}$ denote the resolvent for the operator $H_0,$ $H$ respectively. We factorize the potential
$V=vI_2=V_1V_2$ as for example in (\ref{Y0}).  Later we will show that we can choose $V_2=V.$ Let $Y_0(\l)=V_2R_0(\l)V_1,$ $Y(\l)=V_2R(\l)V_1.$ 

Then we have
%\begin{equation}\lb{res-id}
%Y(\l)=Y_0(\l)-Y_0(\l)\left[I+Y_0(\l)\right]^{-1}Y_0(\l),\qq Y=I-(1+Y_0)^{-1}
%\end{equation} 
%and
\begin{equation}\lb{2.5}
 (I+Y_0(\l))(I-Y(\l))=I.
\end{equation}

As $Y_0(\l):=V_2R_0(\l)V_1\in \cB_2$ is Hilbert-Schmidt but is not trace class (see \cite{Korotyaev2014}),
 we define the modified Fredholm determinant
$$
D(\l)=\det\left[ (I+Y_0(\l)) e^{-Y_0(\l)}\right],\qq\l\in\C_+.
$$

The proofs of the following Corollary and Lemma are identical with the regular case $\vk=0$ and can be found in \cite{IantchenkoKorotyaev2013a}.
\begin{corollary}
\lb{Ttrace}
Let $V\in L^2(\R_+)$ and let $\Im\l\neq 0$. Let $C_\lambda$ be as in (\ref{CL}). Then\\
\no i)
\[
\lb{tr1}\|VR_0(\l)\|_{\cB_2}^2\leq C_\lambda \|V\|_2^2,
\]
\no ii) The operator $R(\l)-R_0(\l)$ is of trace class and satisfies
\[
\lb{tr2}\|R(\l)-R_0(\l)\|_{\cB_1}\leq C_\lambda.
\]

iii) Let, in addition, $V=V_1V_2\in L^2(\R_+)$ with $V_1,V_2\in L^2(\R_+).$ Then for each $\epsilon >0,$
 we have $Y_0, Y, Y_0', Y'\in AC(\cB_2),$  and the following estimates are
 satisfied:
\[
\lb{eY0} \|Y_0(\l)\|_{\cB_2}\le \frac{c}{\epsilon}\|V_1\|_2\|V_2\|_2,\qq \forall \ \l\in \cZ_\epsilon,
\]
\[
\lb{eY0'} \|Y_0(\l)\|_{\cB_2}+\|Y_0'(\l)\|_{\cB_2}\rightarrow 0\qq\mbox{as}\,\,|\Im\l|\rightarrow\infty.
\]

\end{corollary}

\begin{lemma}
\lb{TD1} Let $V\in L^2(\R).$ Then the following facts hold true.

i) For each $\epsilon>0,$ the function $D$ belongs to $AC(\C)$ and satisfies:
\[
\lb{ED1}
D'(\l)=-D(\l)\Tr\left[Y(\l)Y_0'(\l)\right]\qq\forall\l\in \C_+;
\]
\[
\lb{ED2} |D(\l)|\le e^{\|Y_0\|_{\cB_2}},\qq\forall\l\in \C_+;
\]
\[\lb{ED10} D(\l)\rightarrow 1\qq\mbox{as}\qq\Im\l\rightarrow\infty.
\]
ii) For each $\epsilon>0,$ the functions $\log D(\l)$ and $\frac{d}{d\l}\log D(\l)$ belong
to $AC(\C),$ and the following identities  hold true:
\[
\begin{aligned}
\lb{SD} &-\log D(\l)=\sum_{n\ge 2}{\Tr (-Y_0(\l))^n\/n},
\end{aligned}
\]
 where the series converges absolutely and uniformly for $\l$ in the domain
 $$\cL=\{\l\in\C;\,\,\Im \l>c\|V\|_{\cB_2}^2\}$$ for some constant $c>0$ large enough,
  and
\[
\lb{SD2}\begin{aligned}& \rt|\log D(\l)+\sum_{n\ge 2}^N{\Tr (-Y_0(\l))^n\/n}\rt|\le
{\ve_\l^{N+1}\/(N+1)(1-\ve_\l)},\,\,\l\in \cL,\qq
 \ve_\l=C_\l^\frac12(\l) \left\|V\right\|_{\cB_2},
\end{aligned}
\]
for any $N\ge 1.$
Here $C_\l$ is given in (\ref{CL}).  Moreover, $\frac{d^k}{d\l^k}\log D(\l)\in AC
(\C)$ for any $k\in\N$ and the function $D$ is independent of factorization of $V=V_1V_2$ in $Y_0=V_2R_0 V_1,$ so we can choose $Y_0=VR_0.$
\end{lemma}
\vspace{3mm}

We will need the following result on the jump of the cut-off free radial resolvent.

\begin{proposition}\lb{Prop-Jump of Res} Suppose $v\in L^1(\R_+)$  and $\l\in\R,$ $\l\neq\pm m.$ Then the function
$$ \Omega(\l)=\frac{1}{2i}\left(\Tr V (R_0(\l+i0)-R_0(\l-i0))\right)$$ satisfies  (\ref{Oid}) for
 $\l\in (-\infty,-m)\cup (m,+\infty)$ 
\[\lb{OO}
\Omega(\l)=\Omega(\vk,\l)=\int_0^\infty v(y)\left( \frac{k(\l)}{\l-m}\left[kyj_\vk(ky)\right]^2 +\frac{k(\l)}{\l+m}\left[kyj_{\vk-1}(ky)\right]^2 \right)dy,
\] and for $\l\in (-m,m),$ $\Omega(\l)=0.$

Moreover, if in addition $v\in L^\infty(0,a),$ for some $a>0,$ then we have 
\[\lb{Ome}\Omega(\l)=\Omega_0 +{\mathcal O}\left(\frac{\ln |k|}{|k|}\right),\qq\mbox{as}\,\,\l\rightarrow\infty,\qq \mbox{where}\qq \Omega_0=\int_0^\infty v(x)dx.\]

\end{proposition}
{\bf Remark.} The convergence of the integral in (\ref{OO}) follows from (\ref{2.13}).

{\bf Proof of Proposition \ref{Prop-Jump of Res}.}
The integral kernel of the free resolvent $R_0(\l):=(H_0 -\l)^{-1}$ is given by
\[\lb{R_vk1} R_0(x,y,\l)=\left(
      \begin{array}{cc}
        \frac{k(\l)}{\l-m}zh^+_\vk(z)\,\z j_\vk(\z) & zh^+_\vk(z)\,\z j_{\vk-1}(\z) \\
        zh^+_{\vk-1}(z)\,\z j_\vk(\z) & \frac{k(\l)}{\l+m}z h^+_{\vk-1}(z)\,\z j_{\vk-1}(\z) \\
      \end{array}
    \right),\,\, z=kx,\,\, \z=ky,\,\, y<x. \] 
 Note that, for $\vk=0$ formula (\ref{R_vk1}) coincides with the one given in \cite{IantchenkoKorotyaev2013a}, 
as  $zh^+_0(z)=e^{ikx},$ $zh^+_{-1}(z)=ie^{ikx},$ $\z j_0(\z)=\sin(ky)$ and $\z j_{-1}(\z)=\cos(ky).$ 

Recall relations (\ref{Hankel}).
Note the properties which follows from \cite{Antosiewicz1965}, page 439, 10.1.34, 10.1.35:
$$-kxh_\vk^+(-kx)=-kxh_\vk^+(kxe^{i\pi})=-ikx\left(e^{i\vk\pi}j_\vk(kx)+i(-1)e^{i\vk\pi}{}y_\vk(kx)\right).$$
Then, as $k(\l-i0)=-k(\l+i0)$ we have
$$k(\l-i0)xh_\vk^+(k(\l-i0)x)=(-1)^\vk k(\l+i0)h_\vk^-(k(\l+i0)x).$$

Let $y<x$ and $\vk=2,4,\ldots$ even. Then $kyj_\vk(ky)$ is odd with respect to $ky$ and $\frac{k(\l)}{\l-m}kyj_\vk(ky)$ is even, $kyj_{\vk-1}(ky)$ is even. Moreover  $$\begin{aligned}&k(\l-i0)xh_\vk^+(k(\l-i0)x)= k(\l+i0)h_\vk^-(k(\l+i0)x),\\
& k(\l-i0)xh_{\vk-1}^+(k(\l-i0)x)=- k(\l+i0)h_{\vk-1}^-(k(\l+i0)x).\end{aligned}$$

Note also that
$$\begin{aligned}
&kxh_\vk^+(kx)- kxh_\vk^-(kx)=i\left(\frac{\pi k x}{2}\right)^\frac12\left(H_{\vk+1/2}^{(1)}(kx)+H_{\vk+1/2}^{(2)}(kx)\right) =i2(\frac12\pi k x)^\frac12 J_{\vk+1/2}(kx)\\
&=i2kxj_\vk(kx).
\end{aligned}$$

   In order to obtain $R_0(\l+i0)-R_0(\l-i0)$ we calculate
$$\begin{aligned}
& \frac{k(\l)}{\l-m}kxh^+_\vk(kx)\,kyj_\vk(ky)-\frac{k(\l)}{\l-m}kxh^-_\vk(kx)\,kyj_\vk(ky) \\
&=i\frac{k(\l)}{\l-m}2(\frac12\pi k x)^\frac12 J_{\vk+1/2}(kx)\,kyj_\vk(ky)=2i\frac{k(\l)}{\l-m}kx j_\vk(kx)\,kyj_\vk(ky),\\ \\
&kxh^+_\vk(kx)\,kyj_{\vk-1}(ky)-kxh^-_\vk(kx)\,kyj_{\vk-1}(ky)\\
&=i2(\frac12\pi k x)^\frac12 J_{\vk+1/2}(kx)\,kyj_{\vk-1}(ky)=2ikxj_\vk(kx)\,kyj_{\vk-1}(ky),\\ \\
&kxh^+_{\vk-1}(kx)\,kyj_\vk(ky)+(-1)kxh^-_{\vk-1}(kx)\,kyj_\vk(ky)\\
&=i2(\frac12\pi k x)^\frac12 J_{\vk-1/2}(kx)\,kyj_\vk(ky)=2ikxj_{\vk-1}(kx)\,kyj_\vk(ky),\\ \\
&\frac{k(\l)}{\l+m}kxh^+_{\vk-1}(kx)\,kyj_{\vk-1}(ky)+(-1)\frac{k(\l)}{\l+m}kxh^-_{\vk-1}(kx)\,kyj_{\vk-1}(ky)\\
&= \frac{k(\l)}{\l+m}i2(\frac12\pi k x)^\frac12 J_{\vk+1/2}(kx)\,kyj_{\vk-1}(ky)=2i \frac{k(\l)}{\l+m}kxj_{\vk-1}(kx)\,kyj_{\vk-1}(ky)
 \end{aligned}$$
 and
 $$\begin{aligned}
 &V(R_0(\l+i0)-R_0(\l-i0))\\
&= 2iV\left(
                              \begin{array}{cc}
                                \frac{k(\l)}{\l-m}kx j_\vk(kx)\,kyj_\vk(ky) & kxj_\vk(kx)\,kyj_{\vk-1}(ky) \\
                                kxj_{\vk-1}(kx)\,kyj_\vk(ky) & \frac{k(\l)}{\l+m}kxj_{\vk-1}(kx)\,kyj_{\vk-1}(ky)\\
                              \end{array}
                            \right)\theta(x-y).\end{aligned}
  $$
We get
$$\begin{aligned}
&\Tr_{x>y} V (R_0(\l+i0)-R_0(\l-i0))\\= &2i\int_0^\infty v(y)\left( \frac{k(\l)}{\l-m}\left[kyj_\vk(ky)\right]^2 +\frac{k(\l)}{\l+m}\left[kyj_{\vk-1}(ky)\right]^2 \right) dy.
 \end{aligned} $$

The calculation of $\Tr_{x<y} V (R_0(\l+i0)-R_0(\l-i0))$ gives the same formula.
Now, if instead of taking $\vk$ even we suppose that $\vk$ is odd, the rule of  changing of sign in each factor $zh^+_j(z),$  $zj_j(z),$ $\j=\vk,\vk-1,$ will change to the opposite one. As the result the formulas for each entry in the matrix-valued function $(R_0(\l+i0)-R_0(\l-i0)$ will not change.

Note that $$ky{}y_{\vk-1}(ky)=ky (-1)^\vk j_{-\vk}(ky),\qq kyj_0(ky)=\sin ky,\qq ky j_{-1}(ky)=\cos ky.$$

\vspace{5mm}

Now, we prove (\ref{Ome}). 
%Note that for $\l\rightarrow\pm\infty$ $$k(\l)=\l+{\mathcal O}(\l^{-1}),\,\,\frac{\l+m}{k}=1+{\mathcal O}(\l^{-1}),\,\,\frac{k}{\l+m}=1+{\mathcal O}(\l^{-1}).$$
In (\ref{OO}) we split the domain of integration into three intervals
$$\Omega=\int_0^\infty=\int_0^{1/k}+\int_{1/k}^a+\int_a^\infty=\Omega_1+\Omega_2+\Omega_3.$$
In the interval $[0,1/k]$ we apply $|kyj_\vk(ky)|\leq C |ky|^{\vk +1}$ and get
$$\Omega_1(\l)\leq C\int_0^{1/k}|v(y)|\, |ky|^{2\vk}dy\leq C\int_0^{1/k}|v(y)|dy={\mathcal O}(1/|k|)$$  as $v\in L^\infty(0,a).$

We consider the interval $[1/k, a].$
we use
formula (1), page 78, in \cite{Erdelyietal1953II}
\[\lb{F1p78}\begin{aligned}zj_\vk(z)=&\sin \left( z-\vk\frac{\pi}{2}\right)\sum_{j=0}^{\leq\frac12\vk}(-1)^j\left(\vk+\frac12, 2j\right)(2z)^{-2j}+\\
&+\cos\left(z-\vk\frac{\pi}{2}\right)\sum_{j=0}^{\leq\frac12\vk-\frac12}(-1)^j\left(\vk+\frac12, 2j+1\right)(2z)^{-2j-1}.\end{aligned}\]
For $|ky|\geq 1$ formula
(\ref{F1p78}) implies
$kyj_\vk(ky)=\sin\left(ky-\vk\frac{\pi}{2}\right) +{\mathcal O}((ky)^{-1})$ and
\[\lb{OOO}\left[kyj_\vk(ky)\right]^2 +\left[kyj_{\vk-1}(ky)\right]^2=1+{\mathcal O}((ky)^{-1}.\]
Then, as  $v\in L^\infty(0,a),$ we get $$\Omega_2-\int_{1/k}^av(y)dy\leq C\int_{1/k}^a|v(y)|\frac{1}{|ky|}dy=\frac{C}{|k|}\int_{1/k}^a \frac{1}{y}dy={\mathcal O}\left(\frac{\ln |k|}{|k|}\right).$$

In the interval $[1/k, a],$  using  (\ref{OOO}), we get
$$\Omega_3-\int_{a}^\infty v(y)dy\leq C\int_{a}^\infty |v(y)|\frac{1}{|ky|}dy\leq \frac{C}{|k|}.$$

Therefore, we get (\ref{Ome}).
%$$\Omega(\l)=\int_0^\gamma v(y)dy+{\mathcal O}\left(\frac{\ln |k|}{|k|}\right).$$
 \hfill\BBox

\vspace{5mm}

{\bf Proof of Theorem \ref{T1}.} The proof is almost identical to the regular case $\vk=0$ given in \cite{IantchenkoKorotyaev2013a} (see also \cite{IsozakiKorotyaev2011}). We repeat it here for the sake of completeness.

 Let $V\in L^1(\R_+)\cap L^2(\R_+).$

i) We will prove that
${\displaystyle D\in AC(\C),\qq \cS(\l)=\frac{D(\l-i 0)}{D(\l+i 0)}\,e^{-2i\Omega(\l)},\qq\forall\l\in\s_{\rm ac}(H_0),\,\,\l\neq\pm m}.$\\ \\
 Let $\l\in\C_+.$  Denote $\cJ_0(\l)=I+Y_0(\l),$
$\cJ(\l)=I-Y(\l).$ Then $\cJ_0(\l)\cJ(\l)=I$ due to (\ref{2.5}). Now, put
$S_0(\l)=\cJ_0(\overline{\l})\cJ(\l).$ Then we have
$$
S_0(\l)=I-\left(
Y_0(\l)-Y_0(\overline{\l})\right)\left(I-Y(\l)\right).
$$
Now, by the
Hilbert identity,
$$
Y_0(\l)-Y_0(\overline{\l})=(\l-\overline{\l})V_2R_0(\l)R_0(\overline{\l})V_1
$$
 is trace class and
$$\det S_0(\l)=\cS(\l),\qq \l\in \sigma_{\rm ac}(H_0).$$

Let $z=i\tau,$ $\tau\in\R_+$ and $\cD(\l)=\det (\cJ_0(\l)\cJ(z)),$
$\l\in\C_+.$ \\
 It is well defined as $\cJ_0(.)\cJ(z)-I\in
AC(\cB_1).$ The function $\cD(\l)$ is entire in $\C_+$ and
$\cD(z)=I.$ We put $$f(\l)=\frac{D(\l)}{D(z)}e^{\Tr
(Y_0(\l)-Y_0(z))},\qq \l\in\C_+,$$ where
$$D(\l)=\det\left[ (I+Y_0(\l)) e^{-Y_0(\l)}\right].$$ We have
  $\cD(\l)=f(\l),$ $\l\in \C_+.$
Now, using that $\cJ_0(\l)\cJ(\l)=I,$ we get
$$
\det S_0(\l)=\det \cJ_0(\overline{\l})\cJ(z)\cdot
\det(\cJ(z)^{-1}\cJ(\l)=\frac{\cD(\overline{\l})}{\cD(\l)}=
\frac{D(\overline{\l})}{D(\l)}e^{\Tr (Y_0(\overline{\l})-Y_0(\l))}.
$$
As by Proposition \ref{Prop-Jump of Res} we have $\Tr(Y_0(\l+i0)-Y_0(\l-i0))=2i\Omega(\l)$ for
$\l\in \sigma_{\rm ac}(H_0),$   then we get
$$ \cS(\l)=\lim_{\epsilon\downarrow 0}
\frac{D(\l-i\epsilon)}{D(\l+i\epsilon)}e^{-2i\Omega(\l)},\qq\l\in \sigma_{\rm ac}(H_0).$$

Now, by Theorem \ref{th-BFW2},
as $|\l|\rightarrow\infty,$ $\l\in\C_+,$
\[\lb{asJfunc}\gf^+(\l)=-i e^{i\int_0^\infty v(t)dt}+o\left(\l^{-\frac12}\right),\]
and
 we get
 $$\cS(\l)=-\frac{\overline{\gf^+(\l+i0)}}{\gf^+(\l+i0)}
 =e^{-2i\int_0^\infty v(t)dt}+o\left(\l^{-\frac12}\right).$$

ii) We write $\gf(\l)=\gf^+(\l)$ and $\gf^0(\l)=\gf^{0,+}(\l).$
We have $$-\frac{\overline{\gf}(\l+i0)}{\gf(\l+i0)}=\frac{D(\l-i0)}{D(\l+i0)}e^{-2i\Omega(\l)},\qq\l\in\s_{\rm ac}(H_0),$$
where, $$ \Omega(\l)=\frac{1}{2i}\Tr V(R_0(\l+i0)-R_0(\l-i0))\in\R.$$ Now, suppose in addition that $v\in L^\infty(0,a),$ for some $a>0.$ Then by Proposition \ref{Prop-Jump of Res}
$$\Omega(\l)=\int_0^\infty v(t)dt +{\mathcal O}\left(\frac{\ln |k|}{|k|}\right),\qq\l\rightarrow\pm\infty,\qq\l\in\s_{\rm ac}(H_0).$$
Let $\l\in\s_{\rm ac}(H_0)\setminus\{\pm m\}$ and write
$$-\frac{\overline{\gf}(\l+i0)}{\gf(\l+i0)}=\frac{\overline{D(\l+i0)e^{i\Omega(\l)}}}{D(\l+i0)e^{i\Omega(\l)}}\qq\Leftrightarrow
\qq\overline{\left(\frac{\gf(\l+i0)}{D(\l+i0)e^{i\Omega(\l)}}\right)}=-\frac{\gf(\l+i0)}{D(\l+i0)e^{i\Omega(\l)}}.$$
Therefore,
$$e^{i2\arg f_1(0,\l)+i\pi}=e^{i2\arg D(\l)}e^{i2\Omega(\l)},\qq \l\in\sigma_{\rm ac}(H_0)\setminus\{\pm m\}.$$

Moreover, using (\ref{asJfunc})
we get
$$\cS(\l)e^{2i\Omega_0}=\frac{\overline{g}(\l+i0)}{g(\l+i0)}=\frac{\overline{D(\l+i0)e^{i(\Omega(\l)-\Omega_0)}}}{D(\l+i0)e^{i(\Omega(\l)-\Omega_0)}},
\qq\l\in\s_{\rm ac}(H_0),$$
where $$\Omega_0=\int_0^\infty v(t)dt,\qq g(z)=\frac{\gf^+(z)}{k_0(z)e^{i\Omega_0}}.$$

Therefore,
$$e^{i2\arg g(\l)}=e^{i2\arg D(\l)}e^{i2(\Omega(\l)-\Omega_0)},\qq \l\in\sigma_{\rm ac}(H_0)\setminus\{\pm m\}.$$

We know the following facts:\\
  1) $g(\cdot),$ $D(\cdot)\in AC(\C), $ i.e. they are analytic functions on $\C_+,$ continuous in $\overline{\C}_+\setminus\{\pm m\}.$\\
  2) ${\displaystyle g(z) \rightarrow 1,\qq D(z)\rightarrow 1,\qq \Im z\rightarrow\infty},\qq$   $\Omega_0=\int_0^\infty v(x)dx.$

Then the functions $\log g (z),$ $\log D(z)$ are uniquely defined on $\C_+$ and $(-\infty,-m),$ $(m,+\infty);$  and continued from above to the gap $(-m,m).$ Thus
$\log g (z),$ $\log D(z)\in AC(\C)$ and
we have
$$2\arg g(\l)= 2\arg D(\l)+ 2(\Omega(\l)-\Omega_0),\qq\l\in \R\setminus\{\pm m\},\qq\mbox{and}\,\, \Omega(\l)=0\,\,\mbox{for}\,\,\l\in (-m,m).$$

By Cauchy formula, for $z\in\overline{\C}_+\setminus\{\pm m\},$
$$\log g(z)=\frac{1}{\pi}\int\frac{\arg g(t)}{t-z}dt= \frac{1}{\pi}\int\frac{\left(\arg D(t)+\Omega(t)-\Omega_0\right)}{t-z}dt=\log D(z) +\frac{1}{\pi}\int_\R\frac{\left(\Omega(t)-\Omega_0\right)}{t-z}dt,$$
where the  first two integrals are understood in the principal value sense and the last integral is well defined due to $\Omega(t)-\Omega_0={\mathcal O}\left(\frac{\ln |t|}{|t|}\right)\in L^2(\R).$
 Thus
we get \er{a=D}. \phantom{.}\hfill\BBox

\begin{proposition}\lb{P-RR} Suppose all conditions of Theorem \ref{T1} are satisfied. Then for any $\l\in\C_+$
\[\lb{tracemass} \Tr (R(\l)-R_0(\l))=\frac{k_0'(\l)}{k_0(\l)}-\frac{\gf'(\l)}{\gf(\l)},
\] where $k_0(\l)=\gf^0(\vk,\l)$ is the ``free'' radial Jost function.
\end{proposition}
{\bf Proof.}
Using (\ref{ED1}), (\ref{2.5})
$R(\l)-R_0(\l)=-R_0 V_1(I+Y_0 (\l))^{-1}V_2R_0(\l)$ and $Y_0'=V_2R_0^2V_1,$
 we get
\[\lb{DD}\begin{aligned}\frac{d}{d\l} \log D(\l)=\frac{D'(\l)}{D(\l)}&=-\Tr Y(\l)Y_0'(\l)=-\Tr [Y_0'(\l)-(I-Y(\l))Y_0'(\l)]\\
&=-\Tr VR_0^2-\Tr (R(\l)-R_0(\l)).\end{aligned}\]
Now, if potential $v$ satisfies the conditions of Theorem \ref{T1}   then
 \[\lb{D1D} \frac{D'(\l)}{D(\l)}=\frac{\gf'(\l)}{\gf(\l)}-\frac{k_0'(\l)}{k_0(\l)}-\frac{1}{\pi}\int\frac{\Omega(t)-\Omega_0}{(t-\l)^2}dt.\]
Recall that $k_0(\l)=\gf^0(\l)$ is the ``free'' radial Jost function. Now applying
Lemma \ref{res2},
 $$ \Tr (VR_0^2(\l))=\frac{1}{\pi}\int_\sigma\frac{\O(s)}{(s-\l)^2}ds,$$ we get
$$\Tr (R(\l)-R_0(\l))=\frac{k_0'(\l)}{k_0(\l)} -\frac{\gf'(\l)}{\gf(\l)}+\frac{1}{\pi}\int\frac{\Omega(t)-\Omega_0}{(t-\l)^2}dt-\frac{1}{\pi}\int_\sigma\frac{\O(s)}{(s-\l)^2}ds.$$
As $\Omega (t)=0 $ for $t\in (-m,m)$ and $\int\frac{1}{(t-\l)^2}dt=0$ for $\l\in\C\setminus\R,$ we
get (\ref{tracemass}).\hfill\qed

\section{Function $\gF.$ } \lb{s-F}
\setcounter{equation}{0}
In this section we prove Theorems \ref{Th-bound-antibound} and \ref{Prop_counting_zeros}.

\subsection{Characterization of states}

Let $\widetilde{\vt},$ $\widetilde{\vp}$ be solutions of (\ref{Dirac_system}) satisfying
$$ (\widetilde{\vt},\,\,\widetilde{\vp}) = (\vt,\,\, \vp) +o(1)\qq\mbox{as}\,\, x\rightarrow +\infty.$$
By  (\ref{vp0def}),  (\ref{vt0def})  the unperturbed fundamental solutions, $\vp(x,\cdot),\vp(x,\cdot)$  are entire for $x\neq 0$  which implies the following lemma.
\begin{lemma} Let the potential $v$ satisfy (\ref{weakas}). Then 
the functions $\widetilde{\vt}(x,\cdot),$ $\widetilde{\vp}(x,\cdot)$  are entire for each $x> 0$.
\end{lemma}

Now, using (\ref{Jost_Weyl})  we get
\[\lb{Jost_Weyl_pert}\begin{aligned}&f^+(x,\l)= k_0(\l)\widetilde{\vt}(x,\l)  + k^{2\vk}  \widetilde{\vp}(x,\l)  = k_0(\l)\left(\widetilde{\vt} +m_\vk(\l)\widetilde{\vp}\right),\\
 &f^-(x,\l)= \overline{ f^+(x,\l)},\,\,\l\in\sigma_{\rm ac}(H_0).\end{aligned}\] We see that  all singularities of $f^+$  coincide with the singularities of $k_0(\l)$ and do not depend on $x>0.$
As $k_0=\frac{\l+m}{i k(\l)}=\frac{\sqrt{\l+m}}{i\sqrt{\l-m}}$ the only such singularity is at $\l=m.$

Now, the integral kernel of the  resolvent $R(\l):=(H -\l)^{-1}$ is given by
$$ R(x,y,\l)=\left\{\begin{array}{lr}
                 \frac{1}{\det(f^+,\phi)}f^+(x,\l)(\phi(y,\vk,\l))^T & \mbox{if}\,\, y<x, \\
                 \frac{1}{\det(f^+,\phi)}\phi(x,\l)(f^+(y,\vk,\l))^T & \mbox{if}\,\, x<y,
               \end{array}\right.
$$
where $\phi(x,\l)$ is solution of (\ref{Dirac_system}) satisfying (\ref{phi0}). We have $$\begin{aligned}\det(f^+,\phi)=&\lim_{x\rightarrow 0}\frac{x^\vk}{(2\vk-1)!!}f_1^+(x,\l)=\gf^+(\l)\\=&k_0(\l)\lim_{x\rightarrow 0}\frac{x^\vk}{(2\vk-1)!!}\widetilde{\vt}_1(x,\l)  + k^{2\vk}  \lim_{x\rightarrow 0}\frac{x^\vk}{(2\vk-1)!!}\widetilde{\vp}_1(x,\l).\end{aligned}$$
As $\Phi$ is entire, the essential part of the resolvent is

$$\begin{aligned} \mR(x,\l)=&\frac{k_0(\l)\widetilde{\vt}(x,\l)  + k^{2\vk}  \widetilde{\vp}(x,\l) )}{k_0(\l)\lim_{x\rightarrow 0}\frac{x^\vk}{(2\vk-1)!!}\widetilde{\vt}_1(x,\l)  + k^{2\vk}  \lim_{x\rightarrow 0}\frac{x^\vk}{(2\vk-1)!!}\widetilde{\vp}_1(x,\l)}\\
=&\frac{\widetilde{\vt}(x,\l)  + i k^{2\vk+1}(\l+m)^{-1}  \widetilde{\vp}(x,\l) )}{\lim_{x\rightarrow 0}\frac{x^\vk}{(2\vk-1)!!}\widetilde{\vt}_1(x,\l)  +   i k^{2\vk+1}(\l+m)^{-1}   \lim_{x\rightarrow 0}\frac{x^\vk}{(2\vk-1)!!}\widetilde{\vp}_1(x,\l)}
.\end{aligned}$$
The singularities of $\mR(x,\l)$ are independent of $x$ and are
either  zeros of the Jost function $\gf^+(\l)$ or
 $\l=\pm m.$

Note that  for $\vk\in\Z_+,$ $ik^{2\vk+1}(\l+m)^{-1} =k(\l)(\l-m)^\vk (\l+m)^{\vk-1}=0$ at $\l=\pm m$ for $\vk\in \Z_+,$ which is different from the regular case $\vk=0$ as in
\cite{IantchenkoKorotyaev2013} . Otherwise, similar to \cite{IantchenkoKorotyaev2013a}, \cite{IantchenkoKorotyaev2011} and \cite{KorotyaevSchmidt2012}), we get the following  equivalent
characterization of $\s(H).$
\begin{lemma}Let $\vk\in\Z_+.$\\ \no 1)  A point $\l_0\in g^+$ is an eigenvalue iff $\gf^+(\l)=0.$ \\
\no 2) A point $\l_0\in\L_1^-$ is a resonance iff $\gf^+(\l)=0.$\\
\no 3) The  multiplicity of an eigenvalue or a resonance is the multiplicity
of the corresponding zero.\\
\no 4) The point $\l_0=m$ or $\l_0=-m$ is a virtual state iff $\lim_{x\rightarrow 0}\frac{x^\vk}{(2\vk-1)!!}\widetilde{\vt}_1(x,\l_0)=0.$

\end{lemma}

Note that in unperturbed case $\frac{x^\vk}{(2\vk-1)!!}\vt_1(x,\l_0)=(2\vk-1)!!\neq 0,$ which follows from (\ref{vt0def}) and (\ref{bessel2at0}).

\subsection{ Properties of function $\gF$}
We start with some notations. For a function $g=g(\l,x),$ $\l\in\C,$ $x\geq 0,$ we denote $\dot{g}=\partial_\l g,$  $g'=\partial_x g$ and  $g^*(\l):=\overline{g(\overline{\l})}.$

Now, as in \cite{IantchenkoKorotyaev2013a} and similar to \cite{IantchenkoKorotyaev2011} and \cite{KorotyaevSchmidt2012} we introduce an entire function whose zeros contain the states of $H.$
We define $$F(x,\l)=(\l-m)f_1^+(x,\l)f_1^-(x,\l),\qq f_1^-(x,\l)=\left(f_1^+(x,\l)\right)^*,$$ and
\[\lb{function-F}
\gF(\l)=\lim_{x\rightarrow 0}\frac{x^{2\vk}}{((2\vk-1)!!)^2}F(x,\l)=(\l-m)\gf^+(\l)\gf^-(\l),\qq \gf^-(\l)=\left(\gf^+(\l)\right)^* .\]

Such a function was successfully used for the perturbed periodic  Schr{\"o}dinger and Jacobi operators with arbitrary number of gaps (see \cite{KorotyaevSchmidt2012}  and\cite{IantchenkoKorotyaev2011}). 

We have
$$F=(\l-m)\left(k_0(\l)\widetilde{\vt}_1 +k^{2\vk}\widetilde{\vp}_1\right)\left(k_0^*(\l)\widetilde{\vt}_1 +k^{2\vk}\widetilde{\vp}_1\right)=(\l+m)\left(\widetilde{\vt}_1+m_\vk(\l)\widetilde{\vp}_1\right) \left(\widetilde{\vt}_1+m^*_\vk(\l)\widetilde{\vp}_1\right).$$
 Using that for $\l\in\s_{\rm ac}(H_0)$ we have $k^*(\l)=k(\l)$ and
$$ k_0(\l)=\frac{\l +m}{i k(\l)},\qq k_0^*(\l)=-k_0(\l),$$
$$m_\vk(\l)\overline{m}_\vk(\l)=(\l^2-m^2)^{2\vk}\frac{\l-m}{\l+m},\qq m_\vk(\l)+\overline{m}_\vk(\l)=0,$$ we get
\[\lb{Fanother}
F(x,\l)=(\l+m)\widetilde{\vt}_1^2 (x,\l) + (\l-m)(\l^2-m^2)^{2\vk}\widetilde{\vp}_1^2 (x,\l)\] and in unperturbed case $H=H_0$ we have $\gF=\gF_0=(\l-m)k_0k_0^*=\l+m.$
Below we summarize the properties of the function ${\gF}$ and its zeros, generalizing the similar results from \cite{IantchenkoKorotyaev2013a}   for the regular case $\vk=0$    to the present (irregular) problem $\vk\in\Z.$
\begin{proposition}\lb{P-F1}
Assume that potential $v$ satisfies Condition A. Then  function $F$ has the following properties:\\
i)  $\gF(\cdot)$
is entire. \\
ii) $\gF(\cdot)$ is real on  $\R.$  The set of  zeros of $\gF$ is symmetric with respect to the real line. Moreover, $\gF(\l)> 0$ for $\l\in ]-\infty, -m[\cup]m,+\infty[,$  and $\gF$ can have only even number of zeros in $[-m,m].$ \\
iii) If $\l_1$ is an eigenvalue of $H$ then \[\lb{positive}
\dot{\gF}(\l_1)=-2|k(\l_1)|\,\frac{d(\l_1)}{\left(\lim_{x\rightarrow 0}\frac{x^\vk}{(2\vk-1)!!}f_2^+(x,\l_1)\right)^2}<0,\] for some positive function $d(\l_1).$
\end{proposition}
{\em Remark.} Note that $$d(\vk,\l_1)=\lim_{x\rightarrow 0}\left(\frac{x^\vk}{(2\vk-1)!!}\right)^2\int_x^\infty((f_1^+(t,\vk,\l_1))^2+(f_2^+(t,\vk,\l_1))^2)dt$$ and
 in the regular case $\vk=0$ as in \cite{IantchenkoKorotyaev2013a} we have $d(0,\l_1)=
\| f^+(\cdot,\l_1)\|^2_{L^2}.$

{\bf Proof.}
Properies i), ii) follow from formula (\ref{Fanother}) and definition of $\gF$ in (\ref{function-F}).

 The proof of iii) is based on the following  result which can be checked by direct calculation:\\
{\em If $f=(f_1,f_2)^{\rm T}=f(x,\l)$ is solution of the Dirac equation $Hf=\l f$ (\ref{Dirac_system}), then
\[\lb{norming constant} \left(\det (\dot{f},f)\right)'=f_1^2+f_2^2\qq \mbox{for any}\,\, x\in(0,+\infty)\,\,\mbox{and}\,\, \l\in \C\setminus\{\pm m\}.\]}

Now, we fix  $f=f^+(x,\l),$ where  $\l\in g^{+}$  (the upper rim of the gap $(-m,m)$ in $\C\setminus[-m,m]$). Then, as $f^+(\cdot,\vk,\l_1)\in L^2(\R_+,\C^2)$ and  the left hand side of  (\ref{norming constant})  $$\left|
    \begin{array}{cc}
      \dot{f}_1(x,\l_1) & f_1(x,\l_1) \\
       \dot{f}_2(x,\l_1)  & f_2(x,\l_1) \\
    \end{array}
  \right|=\det(\dot{f}(x,\l),f(x,\l))\rightarrow 0\qq\mbox{as}\,\, x\rightarrow\infty,\qq\l\in g^{+},$$  then we get   \[\lb{nc}\det(\dot{f}(x,\l_1),f(x,\l_1))=-\int_x^\infty(f_1^2(t,\l_1)+f_2^2(t,\l_1))dt.\]

 Now, let $\l=\l_1\in\sigma_{\rm bs}(H)$ be an eigenvalue. By applying   $\lim_{x\rightarrow 0}\left(\frac{x^\vk}{(2\vk-1)!!}\right)^2 $ to the both sides of (\ref{nc}), we get  \[\lb{nc2}\begin{aligned}&\dot{\gf}^+(\l_1)\lim_{x\rightarrow 0}\frac{x^\vk}{(2\vk-1)!!}f_2^+(x,\l_1)\\
&=-\lim_{x\rightarrow 0}\left(\frac{x^\vk}{(2\vk-1)!!}\right)^2\int_x^\infty((f_1^+(t,\l_1))^2+(f_2^+(t,\l_1))^2)dt=:-d(\l_1).\end{aligned}\]

Now, going back to the definition of function $\gF$ and differentiating with respect to $\l,$ we get also
$\dot{\gF}(\l_1)=(\l_1 -m)\dot{\gf}^+(\l_1)\gf^-(\l_1).$   Using the Wronskian identity (\ref{wrid}) we get
$$-\lim_{x\rightarrow 0}\frac{x^\vk}{(2\vk-1)!!}f_2^+(x,\l_1)\gf^-(\l_1)=2k_0(\l_1)(k(\l_1))^{2\vk}$$ and therefore
\[\lb{121}
\dot{\gF}(\l_1)=(\l_1 -m)\dot{\gf}^+(\l_1)\frac{-2k_0(\l_1)(k(\l_1))^{2\vk}}{\lim_{x\rightarrow 0}\frac{x^\vk}{(2\vk-1)!!}f_2^+(x,\l_1)}f_2^+(\l_1),\qq k_0(\l_1)=-\frac{\l_1+m}{\sqrt{m^2-\l_1^2}}.\]
Now, putting in (\ref{121}) the expression of   $\dot{\gf}^+(\l_1)$ from equation (\ref{nc2}) we get (\ref{positive}).\hfill\BBox

{\bf Proof of Theorem \ref{Th-bound-antibound}.} 2) follows from the Wronskian identity (\ref{wrid})
 which implies that if $f_1^+(\l_1)=0,$ $\l\neq -m,$ then $f_1^-(\l_1)\neq 0.$

3) follows from identity (\ref{positive}), Proposition \ref{P-F1}.
 \hfill\BBox

 In Proposition \ref{P-F1}  we showed that $F$ is entire in $\C.$ Now, Theorem \ref{Th-ass2terms}, ii), implies that $F$ is
 of exponential type $2\gamma.$

We recall that a function $f$ is
said to  belong    to the Cartwright class $\textsl{Cart}_{\rho_+,\rho_-}$ if $f$ is entire, of
exponential type, and satisfies
$$
\rho_\pm=\rho_\pm(f)\equiv\lim\sup_{y\rightarrow\infty}\frac{\log |f(\pm iy)|}{y} >0,\qq \int_\R\frac{\log(1+|f(x)|)}{1+x^2}dx
<\infty.$$

We determine the asymptotics of the counting function.
We denote $\cN (r,f)$ the total number of zeros of $f$ of modulus $\leq r$ (each zero being counted according to its multiplicity).

We also denote $\cN_+
(r,f)$ (or $\cN_-
(r,f)$) the
number of zeros of function $f$  counted in $\cN (r,f)$ with non-negative (negative) imaginary part  having modulus  $\leq r,$  each zero being counted according to its multiplicity.

\begin{proposition}
\lb{Prop_counting_zeros_F}  Assume that potential $V$ satisfies Condition A and $V'\in L^1(\R_+).$ Then $\gF\in\textsl{Cart}_{2\g,2\g}.$ The set of  zeros of $\gF$ is symmetric with respect to the real line. The set of zeros of $\gF$ with negative imaginary part (i.e. the set of resonances)  satisfies:
\[
\lb{counting}
\cN(r,\gF)= 2 \cN_-(r,\gf^+(\l))={4r\g\/ \pi }(1+o(1))\qq as \qq r\to\iy.
\]
For each $\d >0$ the number of zeros of $\gF$ with negative imaginary part with modulus $\leq r$
lying outside both of the two sectors $|\arg z |<\d ,$ $|\arg z -\pi
|<\d$ is $o(r)$ for large $r$.

\end{proposition}

\section{Massless case. } \lb{s-massless}
\setcounter{equation}{0}

In this section we consider the special case $m=0.$ Then $k=\l,$ $k_0=-i,$ $\s_{\rm ac}(H_0)=\R,$ and the Riemann surface consists of two disjoint sheets $\C$ (see \cite{IantchenkoKorotyaev2013}). This implies that the Jost functions are analytic on $\C$ (see \ref{LemmaH}), and therefore there is no need to introduce a new function $F$ as in the previous section. In the massless case  more results are available.

Note that the massless radial Dirac operator was already studied in  \cite{Blancarteetal1995} and we recall these results. 

\begin{lemma}\lb{L-BGW} Suppose $v\in L^1(\R_+)$ and $m=0.$ Then

\no 1)
the only possible zero of the Jost function $\gf^+(\cdot)$ is $\l=0;$\\
\no 2) $\gf^+(0)=0$ if and only if $\l=0$ is an eigenvalue of $H.$

Moreover, if potential $v$ satisfies (\ref{weakas} ) and $\gf^+(0)=0,$ then
$$\gf^+(\l)=d_\vk\l+{\mathcal O}(|\l|^2)\qq\mbox{as}\,\,|\l|\rightarrow 0,$$ where
$$d_\vk=\frac{i}{c_\vk}\| (\phi(\cdot,\vk,0))\|^2_{L^2}\neq 0,\qq c_\vk=\frac{1}{(2\vk-1)!!}\int_0^\infty x^\vk\phi_2(x,0)v(x)dx,$$
and  $\phi(x,x)$ is regular solution of $H,$ respectively $H_0,$ defined in (\ref{phi0}), respectively (\ref{phi0}).

\end{lemma}
Now,  Corollary \ref{Cor-jostsolutions}, Theorem  \ref{Th-ass2terms}, Lemma \ref{L-BGW}  and Hadamard factorization (see Section 2, equation (2.1) in \cite{IantchenkoKorotyaev2013}) imply
\begin{lemma}\lb{LemmaH} Assume that potential $v$ satisfies Condition A and $v'\in L^1(\R_+).$

The Jost functions $\gf^\pm (\l)$ are entire on $\C.$
Moreover, $\gf^+ (\cdot)\in\textsl{Cart}_{0,2\g},$ and
\[\lb{Had-fact}
\gf^+ (\l)=\l^\sigma c_\vk e^{i\g\,\l }\lim_{r\to +\infty} \prod_{|z_n|\le r}
\lt(1-{\l \/ \l _n} \rt), \ \ \ \l \in \C,\qq \sigma\in\{0,1\} ,
\] where
the product converges uniformly in every bounded disc and
\begin{equation}
\lb{sumcond}
\sum_{z_n\neq 0} {|\Im z_n |\/|z_n|^2} <\infty.
\end{equation}

Here, $\sigma=0$ and $ c_\vk=\gf^+ (0)$ if  $\gf^+ (0)\neq 0,$ and $\sigma=1$ otherwise.   \end{lemma}

We suppose that $v$ satisfies Condition A.
Recall that from Corollary \ref{Ttrace}, (\ref{tr2}), it follows that $R(\l)-R_0(\l)$ is trace class.
Therefore
 $f(H)-f(H_0)$ is trace class for any $f\in {\mathscr S},$ where $\mS$ is the Schwartz class of all rapidly decreasing functions, and
the Krein's trace formula is valid (general result):
$$
\Tr (f(H)-f(H_0))=\int_\R\xi(\l)f'(\l)d\l,\qq f\in \mS,$$ where
$\xi(\l)=\frac{1}{\pi}\phi_{\rm sc}(\l)$ is the spectral shift function and $\phi_{\rm sc}(\l)=\arg \gf^+(\l )+\pi/2=\frac{i}{2}\log\cS$  is the scattering phase.

Let $\gf(\l)=\gf^+(\l)$ be the Jost function. As for $\l\in\R$ the scattering matrix $\cS(\l)$   is given by
$$\cS(\l)=-\frac{\overline{\gf(\l+i0)}}{\gf(\l+i0)}=e^{-2i\phi_{\rm sc}},$$
 then we have also
$$\frac{\cS'}{ \cS}= -2i\Im\frac{\gf'(\l)}{\gf(\l)}.
$$

By using the Hadamard factorization (\ref{Had-fact}) from Lemma  \ref{LemmaH} we get
\[\lb{had}
\frac{\gf'(\l)}{\gf(\l)}=i\g+\lim_{r\to +\infty} \sum_{|\l_n|\le
r}\frac{1}{\l-\l_n},\]

Now, using (\ref{had}) , we get
$$\begin{aligned}&\Tr (f(H)-f(H_0))=\frac{1}{2\pi i}\int_\R f(\l)\frac{\cS'}{ \cS}d\l\\&=-\frac{\gamma}{\pi}\int_{\R}f(\lambda)d\lambda-\frac{1}{\pi}\lim_{r\to +\infty} \sum_{|\l_n|\le
r}\int_\R f(\l)\Im \frac{1}{\l-\l_n}d\l\end{aligned}$$ and
$$\Tr (f(H)-f(H_0))=-\frac{\gamma}{\pi}\int_{\R}f(\lambda)d\lambda-\frac{1}{\pi}\lim_{r\to +\infty} \sum_{|\l_n|\le
r}\int_\R f(\l)\frac{\Im\l_n}{|\l-\l_n|^2}d\l,
$$ 
recovering the Breit-Wigner profile ${\displaystyle -\frac{1}{\pi}\frac{\Im\l_n}{|\l-\l_n|^2}.}$ The sum is converging absolutely by (\ref{sumcond}).

Now applying Proposition \ref{P-RR} with $k_0=-i,$
$$ 
\Tr (R(\l)-R_0(\l))=-\frac{\gf'(\l)}{\gf(\l)},
$$
 and the Hadamard factorization (\ref{had}) in (\ref{DD})
we get the trace formula
$$ \Tr (R(\l)-R_0(\l))=-i\g-\lim_{r\to +\infty} \sum_{|\l_n|\le
r}\frac{1}{\l-\l_n}$$ with uniform convergence in every disc or bounded subset of the plane.

Therefore, the formulas (\ref{TrF1}), (\ref{sc_phase}) and (\ref{TrF3}) in Theorem \ref{T4} are proven.

%\newpage

%\begin{thebibliography}{9}
%\footnotesize

%\end{thebibliography}

%\newpage
\section{Appendix, Proof of Lemma \ref{L-2}}\lb{S-proofL-2}
%\section{Proof of Lemma \ref{L-2}}\lb{S-proofL-2}

We consider equation
\[\lb{inteq}\begin{aligned}
Y(x)=&e^{ikx}\ma1+A \\0 \am+ e^{ikx-i2v_0}\ma
0 \\ B \am        
+\int_x^\g e^{ik\sigma_3(x-t)}(\mN(t)+\mM(t))Y(t)dt.
\end{aligned}
\]
where $v_0=\int_0^\g v(x)dx$. We have 
$$
\mN=i\frac{\l-k}{k}v(t)\sigma_3={\mathcal O}(\l^{-2}),\qq M(t,\l)=e^{-i2\int_0^t v(x)dx}\rt({\vk\/t}-{im\/k} v(t)\rt).
$$ 
Note that
\[\lb{commutation}
e^{-ikt\sigma_3}\mN=\mN e^{-ikt\sigma_3},\qq e^{ikt\sigma_3}\mM=\mM e^{-ikt\sigma_3}.
\]
In the last term in (\ref{inteq}) we use the second commutation relation in (\ref{commutation}) and integrate by parts:
$$
\begin{aligned}
Z=&\int_x^\gamma e^{i\sigma_3 k(x-t)}\mM(t)Y(t)dt=\int_x^\gamma e^{i\sigma_3k(x-2t)}\mM(t)\left(e^{-i\sigma_3 kt}Y(t)\right)dt\\
=&\left[ -\frac{1}{2ik}\sigma_3e^{i\sigma_3k(x-2t)}\mM(t)\left(e^{-i\sigma_3 kt}Y(t)\right)\right]_{t=x}^\gamma+\\&\frac{1}{2ik}\int_x^\gamma
e^{i\sigma_3k(x-2t)}\left\{\mM'(t)e^{-i\sigma_3 kt}-\mM e^{-ikt\sigma_3}(\mN+\mM)\right\}Y(t)dt,
\end{aligned},
$$
where we used that $\widetilde{X}=e^{-ikt\sigma_3}Y$ satisfies
$$
\widetilde{X}'=-\widetilde{W}\widetilde{X},\qq \widetilde{W}=e^{-ikt\sigma_3} (\mN+\mM)e^{ikt\sigma_3}. 
$$

We have $$\mM'=\left( \begin{array}{cc}
0 & \overline{M}'\\
 M' & 0  \\
 \end{array}
 \right),\qq M'(t)=-e^{-i2\int_0^t v(s)ds}\left[i2v(t)\frac{\vk}{t}+2\frac{m}{k} v^2(t)+\frac{\vk}{t^2}+\frac{im}{k} v'(t)\right],$$
\[\lb{c_1} |M'(t)|\leq\frac{c_1}{t^2},\qq 
c_1=\sup_{t\geq 0}\left(t^2M'(t)\right).
\]

By using (\ref{commutation}) and $\mM^2=|M|^2I_2$ we get
$$
Z= \frac{1}{2ik}\sigma_3\mM(x)Y(x)+
\frac{1}{2ik}\int_x^\gamma
e^{i\sigma_3k(x-t)}\left(\mM'(t)-\mM\mN-|M|^2\right)Y(t)dt .
$$ 
Substituting it in (\ref{inteq}) we get
$$
\begin{aligned}
&Y(x)=e^{ikx}\left(\begin{array}{c}            1+A \\
0 \\ \end{array}\right)+ e^{ikx-i2\int_0^\gamma v(t)dt}\left(                                                                                             \begin{array}{c}                                                                                               0 \\ B \\                                                                                            \end{array}                                                                                          \right)+\frac{1}{2ik}\sigma_3\mM(x)Y(x)+\\ 
&+\frac{1}{2ik}\int_x^\gamma
e^{i\sigma_3k(x-t)}\left(2ik\mN(t)+\mM'(t)-\mM\mN-|M|^2\right)Y(t)dt .
\end{aligned}
$$
We have 
$$
|M|^2=\frac{\vk^2}{t^2}+\frac{m^2}{|k|^2}v^2(t)-\frac12 mv(t)\frac{\Im k}{|k|^2}.
$$

We denote  
$$
\mW(t)=2ik\mN(t)+\mM'(t)-\mM\mN-|M|^2,\qq \mA(x)=I-\frac{1}{2ik}\sigma_3\mM(x).
$$ Then
$$\mW(t)=\left(
           \begin{array}{cc}
             \mW_{11} & \mW_{12} \\
             \mW_{21} & \mW_{22} \\
           \end{array}
         \right),$$ where
$$\begin{aligned}&\mW_{11}=-2(\l-k)v(t)-\frac{\vk^2}{t^2}-\frac{m^2}{|k|^2}v^2(t)+\frac12 mv(t)\frac{\Im k}{|k|^2},\\
& \mW_{22}=(\l-k)v(t)-\frac{\vk^2}{t^2}-\frac{m^2}{|k|^2}v^2(t)+\frac12 mv(t)\frac{\Im k}{|k|^2},\qq \mW_{12}=\overline{\mW_{21}},\end{aligned}$$
$$\begin{aligned}&\mW_{21}=\\&-e^{-i2\int_0^t v(s)ds}
\left[i2v(t)\frac{\vk}{t}+2\frac{m}{k} v^2(t)+\frac{\vk}{t^2}+\frac{im}{k} v'(t)+\frac{i(\l-m)\vk}{kt}v(t)+\frac{(\l-k)m}{k^2}v^2(t)\right]
.\end{aligned}$$
Therefore
$$\mW_{11}\sim \mW_{22}\sim -\frac{\vk^2}{t^2}+{\mathcal O}(1)\l^{-1}+{\mathcal O}(1)\l^{-2}, $$
$$\mW_{21}=-e^{-i2\int_0^t v(s)ds}\left[\frac{\vk}{t^2}+i2v(t)\frac{\vk}{t}+{\mathcal O}(1)\frac{1}{kt}+{\mathcal O}(1)\l^{-1}+{\mathcal O}(1)\l^{-3}\right]$$
Then $Y$ satisfies
$$\mA (x,k)Y(x)= Y0+\frac{1}{2ik}\int_x^\gamma
e^{i\sigma_3k(x-t)}\mW(t)Y(t)dt,$$
$$ Y^0(x)=e^{ikx}\left(
                      \begin{array}{c}
                        1+A \\
                        0 \\
                      \end{array}
                    \right)+ e^{ikx-i2\int_0^\gamma v(t)dt}\left(
                                                                                             \begin{array}{c}
                                                                                               0 \\
                                                                                               B \\
                                                                                             \end{array}
                                                                                           \right)=e^{ikx}\left(
                                                                                             \begin{array}{c}
                                                                                               1+A \\
                                                                                               e^{-i2\int_0^\gamma v(t)dt}B \\
                                                                                             \end{array}
                                                                                           \right).$$
We have $$\frac{1}{2ik}\sigma_3\mM(x)=\frac{1}{2ik}\left(
                                        \begin{array}{cc}
                                          0 & \overline{M} \\
                                          -M & 0 \\
                                        \end{array}
                                      \right),$$ where  $$|M(x,\l)|=\left| e^{-i2\int_0^x v(s)ds}\left[\frac{\vk}{x}-\frac{im}{k} v(x)\right]\right|\leq\frac{c_0}{x},\qq c_0=\sup_{x\geq 0}|x M(x)|=\sup_{x\geq 0}\left(\vk -\frac{im}{k} xv(x)\right).$$
                                   We have
 $$
 \left|\frac{1}{2ik}\sigma_3\mM(x)\right| <\frac12\,\,\Leftrightarrow\,\, |kx|>c_0,
 $$
Define  
$$
\mA (x)=I-\frac{1}{2ik}\sigma_3\mM(x),
$$ 
and $\sigma_3\mW=-\mW\sigma_3,$ $\sigma_3^2=\sigma_0,$ $\mM^2=|M|\sigma_0,$ we get
$$
\mA^{-1}=\frac{4k^2}{4k^2-|M|^2}\left(I+\frac{1}{2ik}\sigma_3\mM(x)\right).
$$
Using that \[\lb{ba}
\sup_{|kx|>c_0}|\mA^{-1}(x,k)|\leq 2,\]  we get the integral equation
$$
Y(x)=Y^0+\frac{1}{2ik}(\mA(x,k))^{-1}KY,\qq
 Y^0(x)=\mA^{-1}(x,k)e^{ikx}\left(                                                                                             \begin{array}{c}                                                                                               1+A \\                                                                                               e^{-i2\int_0^\gamma v(t)dt}B \\                                                                                             \end{array}                                                                                           \right),$$ $$ KY=\int_x^\gamma
e^{i\sigma_3k(x-t)}\mW(t)Y(t)dt,
$$
where for $|x|\geq\delta,$ $|\mW(t)|\leq c_1t^{-2}\leq c_1\delta^{-2} $
 by (\ref{c_1}).
   By iterating we get
$$
Y=Y^0+\sum_{n\geq 1}Y^n,\qq  Y^n=\frac{1}{(2ik)^n}(\mA^{-1}K)^nY^0.
$$
Let $ t=(t_j)_1^{n}\in \R^{n}$ and $\mD_t(n)=\{x=t_0<t_1< t_2<...< t_{n}<\gamma\}$.
$$
Y^n=\frac{1}{(2ik)^n}\int_{\mD_t(n)}\prod_{j=1}^n (\mA(t_{j-1}))^{-1}e^{ik \sigma_3 (t_{j-1}-t_j)}\mW(t_j)(\mA(t_n))^{-1}e^{ikt_n}\left(                                                                                             \begin{array}{c}                                                                                               1+A(kt_n) \\                                                                                               e^{-i2\int_0^\gamma v(s)ds}B(kt_n) \\                                                                                             \end{array}                                                                                           \right)dt.
$$                                                                                           Put $\Omega_{\epsilon,\delta}=\{(k,x)\in \cZ_\epsilon^+\times\R_+;\qq \min\{|k|x, x\}>\delta\}$                                                                                           Now, usi(\ref{ba}) we get
$$
 |Y^n(x,k)|\leq \frac{2C_{\epsilon,\delta}}{|k|^n}e^{|\Im k|(2\gamma-x)}\int_{\mD_t(n)}\prod_{j=1}^n |\mW(t_j)|dt=\frac{2}{n!|k|^n}e^{|\Im k|(2\gamma-x)}\left(\int_0^\gamma |\mW(s)|ds\right)^n,  
$$
 where
$$
  C_{\epsilon,\delta}=\sup_{(k,x)\in \C_{\epsilon,\delta}}\{|1+A(kx)|,|B(kx)|\}.
$$\hfill\BBox

 Note that explicitly $$\begin{aligned}&\mA^{-1}=-\frac{2ki}{4k^2-|M|^2}
 \left(  \begin{array}{cc} 2ki & \overline{M} \\
 -M & 2ki \\    \end{array}                                          \right)=b_0\left(                                                        \begin{array}{cc}1 & \frac{1}{2ki}\overline{M} \\   -\frac{1}{2ki}M & 1 \\
 \end{array}
 \right)=\\
&\frac{4k^2}{4k^2-\frac{\vk^2}{t^2}-\frac{m^2}{|k|^2}v^2(t)+\frac12 mv(t)\frac{\Im k}{|k|^2}} \left(                                                                                                \begin{array}{cc}                                                                                                  1 & \frac{1}{2ki}e^{i2\int_0^t v(s)ds}\left[\frac{\vk}{t}+\frac{im}{\overline{k}} v(t)\right] \\                                                                                                  -\frac{1}{2ki}e^{-i2\int_0^t v(s)ds}\left[\frac{\vk}{t}-\frac{im}{k} v(t)\right] & 1 \\                                                                                                \end{array}                                                                                              \right)                                         . \end{aligned}
$$
We write $$\mA^{-1}(x)=b_0\left(
\begin{array}{cc}
 1 & \frac{\vk}{2kxi}e^{i2\int_0^x v(s)ds}\left(1+{\mathcal O}(k^{-1})\right) \\
 -\frac{\vk}{2kxi}e^{-i2\int_0^x v(s)ds}\left(1+{\mathcal O}(k^{-1})\right) & 1 \\
 \end{array}
 \right),
 $$
$$
b_0=\frac{1}{1-\frac14\frac{\vk^2}{(kx)^2}-
\frac{1}{4k^2}\left(\frac{m^2}{|k|^2}v^2(x)-\frac12 mv(x)\frac{\Im k}{|k|^2}\right)}.
$$
If $|k|x\rightarrow\infty$ then
$$b_0= 1+\frac14\frac{\vk^2}{(kx)^2}+{\mathcal O}\left(k^{-3}\right).$$

\setlength{\itemsep}{-\parskip} \footnotesize
\no
{\bf Acknowledgments.} {Various parts of this paper were written
during Evgeny Korotyaev's stay in Aarhus
University, Denmark. He is grateful to the institute for the
hospitality. His study was partly supported by the RFFI grant  No 11-01-00458 and
by  project  SPbGU  No 11.38.215.2014.
}

%\bibliography{coulomb}

\end{document}